# The cubic moment of central values of automorphic $L$-functions

By J. B. Conrey and H. Iwaniec*

**Contents**

## 1. Introduction

The values of $L$-functions at the central point $s = \frac{1}{2}$ (we normalize so that the functional equation connects values at $s$ and $1 - s$) are the subject of intensive studies in various aspects: the algebraicity, the nonvanishing, the positivity. In some instances these numbers are expressible in terms of important geometric invariants (cf. the Birch and Swinnerton-Dyer conjecture for elliptic curves), and the nonvanishing is meaningful in certain structures (such as in the Phillips-Sarnak theory of spectral deformations). The positivity of the Dirichlet $L$-functions with real characters at $s = \frac{1}{2}$ would yield quite remarkable effective lower bounds for the class number of imaginary quadratic fields. Moreover, a good positive lower bound for the central values of Hecke $L$-functions would rule out the existence of the Landau-Siegel zero.

---

*Research of both authors supported by the American Institute of Mathematics and by NSF grants DMS-95-00857, DMS-98-01642.

Independently there is a great interest in upper bounds for the central values; in particular one desires to have a strong estimate in terms of the conductor. The Riemann hypothesis yields the best possible results for individual values, but still there are known unconditional estimates for the average value over distinct families which are as good as the Riemann hypothesis can do, or even slightly better (asymptotic formulas for power moments).

In this paper we consider two families of automorphic $L$-functions associated with the classical (holomorphic) cusp forms of weight $k \geqslant 12$ and the Maass (real-analytic) forms of weight $k = 0$, both for the group $\Gamma = \Gamma_0(q)$ (see the reviews in §2 and §3 respectively). Let $\chi = \chi_q$ be the real, primitive character of modulus $q > 1$. Throughout this paper we assume (for technical simplification) that $q$ is odd, so $q$ is squarefree and $\chi(n) = \left(\frac{n}{q}\right)$ is the Jacobi symbol. To any primitive cusp form $f$ of level dividing $q$ we introduce the $L$-function

$$(1.1) \qquad L_f(s, \chi) = \sum_1^\infty \lambda_f(n)\chi(n)n^{-s}.$$

The main object of our pursuit is the cubic moment

$$(1.2) \qquad \sum_{f \in \mathcal{F}^\star} L_f^3(\tfrac{1}{2}, \chi)$$

where $\mathcal{F}^\star$ is the set of all primitive cusp forms of weight $k$ and level dividing $q$. For this we establish the following bound:

THEOREM 1.1. *Let $k$ be an even number $\geqslant 12$ such that $\chi(-1) = i^k$. Then*

$$(1.3) \qquad \sum_{f \in \mathcal{F}^\star} L_f^3(\tfrac{1}{2}, \chi) \ll q^{1+\varepsilon}$$

*for any $\varepsilon > 0$, the implied constant depending on $\varepsilon$ and $k$.*

Note (see [ILS]) that

$$|\mathcal{F}^\star| = \frac{k-1}{12}\phi(q) + O((kq)^{\frac{2}{3}}).$$

Any cusp form

$$f(z) = \sum_1^\infty \lambda_f(n) n^{\frac{k-1}{2}} e(nz) \in S_k(\Gamma_0(q))$$

yields the twisted cusp form

$$f_\chi(z) = \sum_1^\infty \lambda_f(n)\chi(n) n^{\frac{k-1}{2}} e(nz) \in S_k(\Gamma_0(q^2)),$$

and our $L_f(s, \chi)$ is the $L$-function attached to $f_\chi(z)$. If $f$ is a Hecke form then $f_\chi$ is primitive (even if $f$ is not itself primitive). However, the twisted forms $f_\chi$ span a relatively small subspace of $S_k(\Gamma_0(q^2))$. In view of the above embedding the cubic moment (1.2) looks like a tiny partial sum of a complete sum over the primitive cusp forms of level $q^2$; nevertheless it is alone a spectrally complete sum with respect to the group $\Gamma_0(q)$. This spectral completeness applies more effectively to the normalized cubic moment $\mathcal{C}_k(q)$ which is introduced in (4.14).

In principle our method works also for $k = 2, 4, 6, 8, 10$, but we skip these cases to avoid technical complications. One can figure out, by examining our arguments, that the implied constant in (1.3) is $c(\varepsilon)k^A$, where $A$ is a large absolute number (possibly $A = 3$). We pay some attention to the dependence of implied constants on spectral parameters at the initial (structural) steps, but not in the later analytic transformations. If $q$ tends to $\infty$ over primes, one should be able to get an asymptotic formula

$$(1.4) \qquad \mathcal{C}_k(q) \sim c_k (\log q)^3$$

with $c_k > 0$, but our attempts to accomplish this failed. On the other hand the difficulties of getting an asymptotic formula for $\mathcal{C}_k(q)$ with composite moduli seem to be quite serious (Lemma 14.1 loses the factor $\tau^2(q)$ which causes troubles when $q$ has many divisors; see also Lemma 13.1).

The parity condition $\chi(-1) = i^k$ in Theorem 1.1 can be dropped because if $\chi(-1) = -i^k$ then all the central values $L_f(\frac{1}{2}, \chi)$ vanish by virtue of the minus sign in the functional equation for $L_f(s, \chi)$.

Although our method works for the cubic moment of $L_f(s, \chi)$ at any fixed point on the critical line we have chosen $s = \frac{1}{2}$ for the property

$$(1.5) \qquad L_f(\tfrac{1}{2}, \chi) \geqslant 0.$$

Of course, this property follows from the Riemann hypothesis; therefore it was considered as a remarkable achievement when J.-L. Waldspurger [Wa] derived (1.5) from his celebrated formula; see also W. Kohnen and D. Zagier [KZ]. Without having the nonnegativity of central values one could hardly motivate the goal of estimating the cubic moment (still we would not hesitate to get an asymptotic formula). As a consequence of (1.5) we derive from (1.3) the following bound for the individual values.

COROLLARY 1.2. *Let $f$ be a primitive cusp form of weight $k \geqslant 12$ and level dividing $q$, and let $\chi(\mathrm{mod}\, q)$ be the primitive real character (the Jacobi symbol). Then*

$$(1.6) \qquad L_f(\tfrac{1}{2}, \chi) \ll q^{\frac{1}{3}+\varepsilon}$$

*for any $\varepsilon > 0$, the implied constant depending on $\varepsilon$ and $k$.*

Let us recall that the convexity bound is $L_f(\frac{1}{2}, \chi) \ll q^{\frac{1}{2}+\varepsilon}$ while the Riemann hypothesis yields $L_f(\frac{1}{2}, \chi) \ll q^\varepsilon$.

An interesting case is that for a Hecke cusp form $f$ of level one and weight $k$ (so $k$ is an even integer $\geq 12$),

$$f(z) = \sum_1^\infty a(n) n^{\frac{k-1}{2}} e(nz) \in S_k(\Gamma_0(1)).$$

This corresponds, by the Shimura map, to a cusp form $g$ of level four and weight $\frac{k+1}{2}$,

$$g(z) = \sum_1^\infty c(n) n^{\frac{k-1}{4}} e(nz) \in S_{\frac{k+1}{2}}(\Gamma_0(4)).$$

We normalize $f$ by requiring $a(1) = 1$, while $g$ is normalized so that

$$\frac{1}{6} \int_{\Gamma_0(4) \backslash \mathbb{H}} |g(z)|^2 y^{\frac{k+1}{2}} d\mu z = 1.$$

Then the formula of Waldspurger, as refined by Kohnen-Zagier (see Theorem 1 of [KZ]), asserts that for $q$ squarefree with $\chi_q(-1) = i^k$,

(1.7) $$c^2(q) = \pi^{-\frac{k}{2}} \Gamma(\tfrac{k}{2}) L_f(\tfrac{1}{2}, \chi_q) \langle f, f \rangle^{-1}$$

where

$$\langle f, f \rangle = \int_{\Gamma_0(q) \backslash \mathbb{H}} |f(z)|^2 y^k d\mu z.$$

By (1.6) and (1.7) we get:

COROLLARY 1.3. *If $q$ is squarefree with $\chi_q(-1) = i^k$ then*

(1.8) $$c(q) \ll q^{\frac{1}{6} + \varepsilon}$$

*where the implied constant depends on $\varepsilon$ and the form $f$.*

This result constitutes a considerable improvement of the estimates given in [I1] and [DFI]. It also improves the most recent estimate by V. A. Bykovsky [By] who proved (1.8) with exponent $3/16$ in place of $1/6$. Actually [DFI] and [By] provide estimates for $L_f(s, \chi)$ at any point on the critical line. To this end (as in many other papers; see the survey article [Fr] by J. Friedlander) the second moment of relevant $L$-functions is considered with an amplifier which is the square of a short Dirichlet polynomial. In such a setting the property (1.5) is not needed yet a sub-convexity bound is achieved by proper choice of the length and the coefficients of the amplifier.

Here is the second instance where the nonnegativity of central values of automorphic $L$-functions plays a crucial role for their estimation (the first case appears in [IS] in the context of the Landau-Siegel zero). By comparison with the former methods one may interpret the cubic moment approach as a kind of amplification of $L_f^2(\tfrac{1}{2}, \chi)$ by the factor $L_f(\tfrac{1}{2}, \chi)$. In this role as a self-amplifier the central value $L_f(\frac{1}{2}, \chi)$ is represented by a Dirichlet polynomial whose length exceeds greatly all of these in previous practice (here it has length about $q$). Hence the question: what makes it possible to handle the present case? There are many different reasons; for instance we emphasize the smoothness of the self-amplification. Consequently, it can be attributed to our special amplifier $L_f(\frac{1}{2}, \chi)$ that at some point the character sums $g(\chi, \psi)$ in two variables over a finite field crop (see (11.10)). From then on our arguments are powered by the Riemann hypothesis for varieties (Deligne's theory, see §13). In fact our arguments penetrate beyond the Riemann hypothesis as we exploit the variation in the angle of the character sum (10.7) when estimating general bilinear forms (11.1) (see the closing remarks of §11).

To reduce the spectral sum (1.2) to the character sums in question we go via Petersson's formula to Kloosterman sums, open the latter and execute the resulting additive character sums in three variables (which come from a smooth partition of $L_f^3(\frac{1}{2}, \chi)$ into Dirichlet polynomials) by Fourier analysis on $\mathbb{R}^3$. One may argue that our computations would be better performed by employing harmonic analysis on $\mathrm{GL}_3(\mathbb{R})$; however, we prefer to use only the classical tools (Poisson's formula) which are commonly familiar. In this connection we feel the demand is growing for practical tables of special functions on higher rank groups to customize them as much as the Bessel functions are on $\mathrm{GL}_2(\mathbb{R})$. Still, there is a revealing advantage to direct computations; see our comments about the factor $e(mm_1m_2/c)$ in (8.32) and (10.1), which presumably would not be visible in the framework of $\mathrm{GL}_3(\mathbb{R})$. This technical issue sheds some light on the position of Bessel functions towards Kloosterman sums.

In this paper we also consider the spectral cubic moment of central values of $L$-functions attached to Maass forms of weight zero. Since the space of such forms is infinite we take only those with bounded spectral parameter; i.e., we consider

$$(1.9) \qquad \sideset{}{^\star}\sum_{|t_j| \leqslant R} L_j^3(\tfrac{1}{2}, \chi) + \int_{-R}^{R} |L(\tfrac{1}{2} + ir, \chi)|^6 \ell(r) \, dr$$

where $\ell(r) = r^2(4 + r^2)^{-1}$. We refer the reader to Sections 3 and 5 to find the terminology. Actually the Maass forms were our primary interest when we started. Here the special attraction lies in the subspace of the continuous spectrum which is spanned by the Eisenstein series $E_\mathfrak{a}(z, \frac{1}{2} + ir)$ (there are $\tau(q)$ distinct Eisenstein series associated with the cusps $\mathfrak{a}$ of $\Gamma_0(q)$). Every Eisenstein series gives us the same $L$-function $L(s - ir, \chi)L(s + ir, \chi)$ (however, with different proportions equal to the width of the cusp; see (3.27)) whose central value is $|L(\frac{1}{2} + ir, \chi)|^2$; hence its cube is the sixth power of the Dirichlet $L$-function.

THEOREM 1.4. *Let $R \geqslant 1$. For any $\varepsilon > 0$,*

$$(1.10) \qquad \sideset{}{^\star}\sum_{|t_j| \leqslant R} L_j^3(\tfrac{1}{2}, \chi) + \int_{-R}^{R} |L(\tfrac{1}{2} + ir, \chi)|^6 \ell(r) \, dr \ll R^A q^{1+\varepsilon}$$

*with some absolute constant $A \geqslant 1$, the implied constant depending on $\varepsilon$.*

Here, as in the case of holomorphic cusp forms, the central values $L_j(\tfrac{1}{2}, \chi)$ are also known to be nonnegative without recourse to the Riemann hypothesis due to Katok-Sarnak [KS] and Guo [Gu]. Note that in the space of continuous spectrum this amounts to $|L(\tfrac{1}{2} + ir, \chi)|^2 \geqslant 0$. However, this property for the central values of cuspidal $L$-functions is quite subtle, and is indispensable in what follows, even for estimating the Dirichlet $L$-functions. First it allows us to derive from (1.10) the corresponding extensions of the estimates (1.6) and (1.8) for the Maass cusp forms. Another observation is that we receive the Dirichlet $L$-functions at any point on the critical line (not just at $s = \tfrac{1}{2}$ as for the cuspidal $L$-functions) by virtue of the integration in the continuous spectrum parameter. Ignoring the contribution of the cuspidal spectrum in (1.10) and applying Hölder's inequality to the remaining integral, one derives

$$(1.11) \qquad \int_{-R}^{R} |L(\tfrac{1}{2} + ir, \chi)| \, dr \ll R^A q^{\frac{1}{6}+\varepsilon}$$

where the implied constant depends on $\varepsilon$ (in this way we relax the peculiar measure $\ell(r)dr$ which vanishes at $r = 0$ to order two). Hence, we have the following result:

COROLLARY 1.5. *Let $\chi$ be a real, nonprincipal character of modulus $q$. Then for any $\varepsilon > 0$ and $s$ with $\operatorname{Re} s = \tfrac{1}{2}$,*

$$(1.12) \qquad L(s, \chi) \ll |s|^A q^{\frac{1}{6}+\varepsilon}$$

*where $A$ is an absolute constant and the implied constant depends on $\varepsilon$.*

It would not be difficult to produce a numerical value of $A$ which is quite large. A hybrid bound which is sharp in both the $s$ aspect and the $q$-aspect simultaneously (not only for the real character) was derived by R. Heath-Brown [H-B] by mixing the van der Corput method of exponential sums and the Burgess method of character sums. In the $q$-aspect alone our bound (1.12) marks the first improvement of the celebrated result of D. Burgess [Bu] with exponent $3/16$ in place of $1/6$. Moreover, our exponent $1/6$ matches the one in the classical bound for the Riemann zeta-function on the line $\operatorname{Re} s = \tfrac{1}{2}$, which can be derived by Weyl's method of estimating exponential sums. Though Weyl's method has been sharpened many times (see the latest achievement of M. N. Huxley [Hu]) any improvement of (1.12) seems to require new ideas (we tried to introduce an extra small amplification to the cubic moments without success). On this occasion let us recall that the aforementioned methods of Weyl and Burgess yield the first boundsbreaking convexity for the $L$-functions on $GL_1$. Since then many refinements and completely new methods were developed for both the $L$-functions on $GL_1$ and the $L$-functions on $GL_2$; see [Fr]. We should also point out that Burgess established nontrivial bounds for character sums of length $N \gg q^{\frac{1}{4}+\epsilon}$, while (1.12) yields nontrivial bounds only if $N \gg q^{\frac{1}{3}+\epsilon}$. In particular, (1.12) does not improve old estimates for the least quadratic nonresidue.

Our main goal (proving Corollaries 1.2 and 1.5) could be accomplished in one space $\mathcal{L}_k(\Gamma_0(q)\backslash\mathbb{H})$ of square-integrable functions $F : \mathbb{H} \to \mathbb{C}$ which transform by

$$(1.13) \qquad F(\gamma z) = \left(\frac{cz+d}{|cz+d|}\right)^k F(z)$$

for all $\gamma \in \Gamma_0(q)$. In this setting the holomorphic cusp forms $f(z)$ of weight $k$ (more precisely the corresponding forms $F(z) = y^{-k/2}f(z)$) lie at the bottom of the spectrum, i.e., in the eigenspace of $\lambda = \frac{k}{2}(1 - \frac{k}{2})$ of the Laplace operator

$$(1.14) \qquad \Delta_k = y^2 \left(\frac{\partial^2}{\partial x^2} + \frac{\partial^2}{\partial y^2}\right) - iky\frac{\partial}{\partial x},$$

while the Eisenstein series still yield the Dirichlet $L$-functions on the critical line. We have chosen to present both cases of holomorphic and real-analytic forms separately to illustrate structural differences until the end of Section 5. From this point on both cases are essentially the same so we restrict our arguments to the holomorphic forms.

*Acknowledgement.* Our work on this paper began and was nearly finished in July 1998 at the American Institute of Mathematics in Palo Alto, California. The second author is grateful to the Institute for the invitation and generous support during his visit. He also wishes to express admiration to John Fry for his unprecedented will to support research in mathematics in America and his deep vision of the AIM. Finally, we thank the referee for careful reading and valuable corrections.

## 2. A review of classical modular forms

Let $q$ be a positive integer. We restrict our considerations to the Hecke congruence group of level $q$ which is

$$\Gamma_0(q) = \left\{\begin{pmatrix} a & b \\ c & d \end{pmatrix} \in \mathrm{SL}_2(\mathbb{Z}) : c \equiv 0 (\mathrm{mod}\, q)\right\};$$

its index in the modular group is

$$\nu(q) = [\Gamma_0(1) : \Gamma_0(q)] = q \prod_{p|q}(1 + \frac{1}{p}).$$

The group $\Gamma = \Gamma_0(q)$ acts on the upper half-plane $\mathbb{H} = \{z = x + iy : y > 0\}$ by

$$\gamma = \frac{az+b}{cz+d}, \quad \text{if } \gamma = \begin{pmatrix} a & b \\ c & d \end{pmatrix} \in \Gamma.$$

Let $k$ be a positive even integer. The space of cusp forms of weight $k$ and level $q$ is denoted $S_k(\Gamma_0(q))$; it is a finite-dimensional Hilbert space with respect to the inner product

$$(2.1) \qquad \langle f, g \rangle = \int_{\Gamma \backslash \mathbb{H}} f(z)\overline{g}(z) y^k d\mu z$$

where $d\mu z = y^{-2} dx dy$ is the invariant measure on $\mathbb{H}$. Let $\mathcal{F} = \{f\}$ be an orthonormal basis of $S_k(\Gamma_0(q))$. We can assume that every $f \in \mathcal{F}$ is an eigenfunction of the Hecke operators

$$(2.2) \qquad (T_n f)(z) = \frac{1}{\sqrt{n}} \sum_{ad=n} (\frac{a}{d})^{k/2} \sum_{b(\bmod d)} f\left(\frac{az+b}{d}\right),$$

for all $n$ with $(n, q) = 1$; i.e., $T_n f = \lambda_f(n) f$ if $(n, q) = 1$. We call $\mathcal{F}$ the Hecke basis of $S_k(\Gamma_0(q))$. The eigenvalues $\lambda_f(n)$ are related to the Fourier coefficients of $f(z)$. We write

$$(2.3) \qquad f(z) = \sum_{1}^{\infty} a_f(n) n^{\frac{k-1}{2}} e(nz).$$

Then for $(n, q) = 1$ we have

$$(2.4) \qquad a_f(n) = a_f(1) \lambda_f(n).$$

Note that if $a_f(1) = 0$ then $a_f(n) = 0$ for all $n$ co-prime with $q$. The Hecke eigenvalues $\lambda_f(n)$ are real and they have the following multiplicative property

$$(2.5) \qquad \lambda_f(m) \lambda_f(n) = \sum_{d|(m,n)} \lambda_f(mn/d^2)$$

if $(mn, q) = 1$.

For any orthonormal basis $\mathcal{F}$ of $S_k(\Gamma_0(q))$ and any $m, n \geq 1$ we have the following Petersson formula (cf. Theorem 3.6 of [I3]):

$$(2.6) \quad (4\pi)^{1-k} \Gamma(k-1) \sum_{f \in \mathcal{F}} \bar{a}_f(m) a_f(n)$$

$$= \delta(m,n) + 2\pi i^k \sum_{c \equiv 0 (\bmod q)} c^{-1} S(m,n;c) J_{k-1}\left(\frac{4\pi}{c}\sqrt{mn}\right)$$

where $\delta(m,n)$ is the Kronecker diagonal symbol, $S(m,n;c)$ is the Kloosterman sum defined by

$$(2.7) \qquad S(m,n;c) = \sum_{ad \equiv 1 (\mathrm{mod}\, c)} e\left(\frac{am+dn}{c}\right),$$

and $J_{k-1}(x)$ is the Bessel function of order $k-1$. Notice that the series on the right-hand side of (2.6) converges absolutely by virtue of the Weil bound

$$(2.8) \qquad |S(m,n;c)| \leq (m,n,c)^{\frac{1}{2}} c^{\frac{1}{2}} \tau(c).$$

For the orthonormal Hecke basis and $(mn,q) = 1$ we can write (2.6) as

$$(2.9)$$
$$\sum_{f \in \mathcal{F}} \omega_f \lambda_f(m) \lambda_f(n) = \delta(m,n) + \sqrt{mn} \sum_{c \equiv 0 (\mathrm{mod}\, q)} c^{-2} S(m,n;c) J(2\sqrt{mn}/c)$$

where

$$(2.10) \qquad \omega_f = (4\pi)^{1-k} \Gamma(k-1) |a_f(1)|^2,$$

$$(2.11) \qquad J(x) = 4\pi i^k x^{-1} J_{k-1}(2\pi x).$$

According to the Atkin-Lehner theory [AL] the sum (2.9) can be arranged into a sum over all primitive forms of level dividing $q$, but, of course, with slightly different coefficients. Precisely, a primitive form $f$ appears with coefficient

$$(2.12) \qquad \omega_f^* = \frac{12}{(k-1)q} \left( \sum_{(\ell,q)=1} \lambda_f(\ell^2) \ell^{-1} \right)^{-1} \gg (kq)^{-1-\varepsilon}$$

for any $\varepsilon > 0$, the implied constant depending only on $\varepsilon$.

*Remarks.* For the formula (2.12) see [ILS]. The coefficient $\omega_f^*$ is essentially (up to a simple constant factor) the inverse of the symmetric square $L$-function associated with $f$ at the point $s = 1$. J. Hoffstein and P. Lockhart [HL] showed that $\omega_f^* \ll (kq)^{\varepsilon-1}$, but we do not need this bound for applications in this paper. The lower bound (2.12) can be established by elementary arguments.

Later we assume that $k \geq 12$ to secure a sufficiently rapid convergence of the series of Kloosterman sums in (2.9). Indeed we have

$$J_{k-1}(x) \ll \min(x^{k-1}, x^{-1/2})$$

which yields

$$(2.13) \qquad J(x) \ll \min(x^{10}, x^{-3/2}) \ll x^{10} (1+x^2)^{-23/4}.$$

## 3. A review of Maass forms

In this section we introduce the notation and basic concepts from the theory of Maass forms of weight $k = 0$ in the context of the Hecke congruence group $\Gamma = \Gamma_0(q)$. There is no essential difference from the theory of classical forms except for the existence of a continuous spectrum in the space of Maass forms. This is important for our applications since it brings us the Dirichlet $L$-function.

Let $\mathcal{A}(\Gamma\backslash\mathbb{H})$ denote the space of automorphic functions of weight zero, i.e., the functions $f : \mathbb{H} \to \mathbb{C}$ which are $\Gamma$-periodic. Let $\mathcal{L}(\Gamma\backslash\mathbb{H})$ denote the subspace of square-integrable functions with respect to the inner product (2.1) with $k = 0$. The Laplace operator

$$\Delta = y^2 \left( \frac{\partial^2}{\partial x^2} + \frac{\partial^2}{\partial y^2} \right)$$

acts in the dense subspace of smooth functions in $\mathcal{L}(\Gamma\backslash\mathbb{H})$ such that $f$ and $\Delta f$ are both bounded; it has a self-adjoint extension which yields the spectral decomposition $\mathcal{L}(\Gamma\backslash\mathbb{H}) = \mathbb{C} \oplus \mathcal{C}(\Gamma\backslash\mathbb{H}) \oplus \mathcal{E}(\Gamma\backslash\mathbb{H})$. Here $\mathbb{C}$ is the space of constant functions, $\mathcal{C}(\Gamma\backslash\mathbb{H})$ is the space of cusp forms and $\mathcal{E}(\Gamma\backslash\mathbb{H})$ is the space of Eisenstein series.

Let $\mathcal{U} = \{u_j : j \geq 1\}$, be an orthonormal basis of the $\mathcal{C}(\Gamma\backslash\mathbb{H})$ which are eigenfunctions of $\Delta$, say

$$(\Delta + \lambda_j)u_j = 0 \quad \text{with} \quad \lambda_j = s_j(1 - s_j), \quad s_j = \frac{1}{2} + it_j.$$

Since $\lambda_j \geq 0$ we have $\operatorname{Re} s_j = \frac{1}{2}$ or $\frac{1}{2} \leq s_j < 1$. Any $u_j(z)$ has the Fourier expansion of type

$$(3.1) \qquad u_j(z) = \sum_{n \neq 0} \rho_j(n) W_{s_j}(nz)$$

where $W_s(z)$ is the Whittaker function given by

$$(3.2) \qquad W_s(z) = 2|y|^{\frac{1}{2}} K_{s-\frac{1}{2}}(2\pi|y|) e(x)$$

and $K_s(y)$ is the $K$-Bessel function. Note that $W_s(z) \sim e(z)$ as $y \to \infty$. The automorphic forms $u_j(z)$ are called Maass cusp forms.

The eigenpacket in $\mathcal{E}(\Gamma\backslash\mathbb{H})$ consists of Eisenstein series $E_\mathfrak{a}(z, s)$ on the line $\operatorname{Re} s = \frac{1}{2}$. These are defined for every cusp $\mathfrak{a}$ by

$$E_\mathfrak{a}(z, s) = \sum_{\gamma \in \Gamma_\mathfrak{a}\backslash\Gamma} (\operatorname{Im} \sigma_\mathfrak{a}^{-1}\gamma z)^s$$

if $\operatorname{Re} s > 1$ and by analytic continuation for all $s \in \mathbb{C}$. Here $\Gamma_\mathfrak{a}$ is the stability group of $\mathfrak{a}$ and $\sigma_\mathfrak{a} \in \operatorname{SL}_2(\mathbb{R})$ is such that $\sigma_\mathfrak{a}\infty = \mathfrak{a}$ and $\sigma_\mathfrak{a}^{-1}\Gamma_\mathfrak{a}\sigma_\mathfrak{a} = \Gamma_\infty$. The scaling matrix $\sigma_{\mathfrak{a}}$ of cusp $\mathfrak{a}$ is only determined up to a translation from the right; however the Eisenstein series does not depend on the choice of $\sigma_{\mathfrak{a}}$, not even on the choice of a cusp in the equivalence class. The Fourier expansion of $E_{\mathfrak{a}}(z, s)$ is similar to that of a cusp form; precisely,

$$(3.3) \qquad E_{\mathfrak{a}}(z,s) = \varphi_{\mathfrak{a}} y^s + \varphi_{\mathfrak{a}}(s) y^{1-s} + \sum_{n \neq 0} \varphi_{\mathfrak{a}}(n,s) W_s(nz)$$

where $\varphi_{\mathfrak{a}} = 1$ if $\mathfrak{a} \sim \infty$ or $\varphi_{\mathfrak{a}} = 0$ otherwise.

We can assume that $\mathcal{U}$ is the Hecke basis, i.e., every $u_j \in \mathcal{U}$ is an eigenfunction of all the Hecke operators (2.2) with $k = 0$,

$$(3.4) \qquad T_n u_j = \lambda_j(n) u_j \quad \text{if} \quad (n,q) = 1.$$

Moreover, the reflection operator $R$ defined by $(Rf)(z) = f(-\bar{z})$ commutes with $\Delta$ and all $T_n$ with $(n,q) = 1$ so that we can also require

$$(3.5) \qquad R u_j = \varepsilon_j u_j.$$

Since $R$ is an involution the space $\mathcal{C}(\Gamma \backslash \mathbb{H})$ is split into even and odd cusp forms according to $\varepsilon_j = 1$ and $\varepsilon_j = -1$. All the Eisenstein series $E_{\mathfrak{a}}(z,s)$ are even and they are also eigenfunctions of the Hecke operators

$$(3.6) \qquad T_n E_{\mathfrak{a}}(z,s) = \eta_{\mathfrak{a}}(n,s) E_{\mathfrak{a}}(z,s), \quad \text{if} \quad (n,q) = 1.$$

The analog of Petersson's formula (2.6) for Maass forms is the following formula of Kuznetsov (see Theorem 9.3 of [I2]):

$$(3.7)$$
$$\sum_j h(t_j) \bar{\rho}_j(m) \rho_j(n) + \sum_{\mathfrak{a}} \frac{1}{4\pi} \int_{-\infty}^{\infty} h(r) \bar{\varphi}_{\mathfrak{a}}(m, \tfrac{1}{2}+ir) \varphi_{\mathfrak{a}}(n, \tfrac{1}{2}+ir) \, dr$$
$$= \delta(m,n) H + \sum_{c \equiv 0 (\bmod q)} c^{-1} S(m,n;c) H^{\pm}\Big(\frac{4\pi}{c}\sqrt{|mn|}\Big)$$

where $\pm$ is the sign of $mn$ and $H, H^+(x), H^-(x)$ are the integral transforms of $h(t)$ given by

$$(3.8) \qquad H = \frac{1}{\pi} \int_{-\infty}^{\infty} h(t) \operatorname{th}(\pi t) t\, dt,$$

$$(3.9) \qquad H^+(x) = 2i \int_{-\infty}^{\infty} J_{2it}(x) \frac{h(t) t}{\operatorname{ch} \pi t} dt,$$

$$(3.10) \qquad H^-(x) = \frac{4}{\pi} \int_{-\infty}^{\infty} K_{2it}(x) \operatorname{sh}(\pi t) h(t) t\, dt.$$

This formula holds for any orthonormal basis $\mathcal{U}$ of cusp forms in $\mathcal{C}(\Gamma \backslash \mathbb{H})$, for any $mn \neq 0$ and any test function $h(t)$ which satisfies the following conditions;

$(3.11) \qquad h(t) \quad \text{is holomorphic in} \quad |\operatorname{Im} t| \leqslant \sigma,$

(3.12) $$h(t) = h(-t),$$
(3.13) $$h(t) \ll (|t|+1)^{-\theta},$$

for some $\sigma > \frac{1}{2}$ and $\theta > 2$.

For the Fourier coefficients $\varphi_{\mathfrak{a}}(\pm n, s)$ of the Eisenstein series $E_{\mathfrak{a}}(z, s)$,

(3.14) $$\varphi_{\mathfrak{a}}(\pm n, s) = \varphi_{\mathfrak{a}}(1, s)\eta_{\mathfrak{a}}(n, s)n^{-\frac{1}{2}},$$

if $n > 0$, $(n, q) = 1$. For the coefficients $\rho_j(n)$ of the Hecke-Maass form $u_j(z)$ we have a similar formula

(3.15) $$\rho_j(\pm n) = \rho_j(\pm 1)\lambda_j(n)n^{-\frac{1}{2}},$$

if $n > 0$, $(n, q) = 1$. Moreover

(3.16) $$\rho_j(-1) = \varepsilon_j \rho_j(1).$$

To simplify presentation we restrict the spectral sum in (3.7) to the even forms; these can be selected by adding (3.7) for $m, n$ to that for $-m, n$. We obtain for $m, n \geqslant 1$, $(mn, q) = 1$,

(3.17)
$$\sideset{}{'}\sum_j h(t_j)\omega_j \lambda_j(m)\lambda_j(n) + \sum_{\mathfrak{a}} \frac{1}{4\pi}\int_{-\infty}^{\infty} h(r)\omega_{\mathfrak{a}}(r)\eta_{\mathfrak{a}}(m, \tfrac{1}{2}+ir)\eta_{\mathfrak{a}}(n, \tfrac{1}{2}+ir)\,dr$$
$$= \tfrac{1}{2}\delta(m, n)H + \sqrt{mn}\sum_{c \equiv 0 (\mathrm{mod}\, q)}$$
$$\times c^{-2}\{S(m, n; c)J^+(2\sqrt{mn}/c) + S(-m, n; c)J^-(2\sqrt{mn}/c)\}$$

where $\sum'$ restricts to the even Hecke cusp forms,

(3.18) $$\omega_j = 4\pi|\rho_j(1)|^2/\operatorname{ch}\pi t_j$$

for $\lambda_j = s_j(1 - s_j)$ with $s_j = \tfrac{1}{2} + it_j$ and

(3.19) $$\omega_{\mathfrak{a}}(r) = 4\pi|\varphi_{\mathfrak{a}}(1, \tfrac{1}{2}+ir)|^2/\operatorname{ch}\pi r.$$

The $J$-functions which are attached to the Kloosterman sums on the right-hand side of (3.17) are defined by $J^{\pm}(x) = x^{-1}H^{\pm}(2\pi x)$. In our applications of (3.17) we assume that the conditions (3.11)-(3.13) hold with $\sigma > 6$ to ensure the bound $H^+(x) \ll \min(x^{11}, x^{-1/2})$. For $x \geqslant 1$ this follows by $J_{2it}(x) \ll x^{-1/2}\operatorname{ch}\pi t$, and for $0 < x < 1$ this follows by moving the integration in (3.9) to the horizontal line $\operatorname{Im} t = 6$ and applying $J_s(x) \ll x^{\sigma}e^{\pi|s|/2}$. The same bound is derived for $H^-(x)$ by similar arguments. In any case we get

(3.20) $$J^{\pm}(x) \ll \min(x^{10}, x^{-3/2}) \ll x^{10}(1 + x^2)^{-23/4}.$$

Recall that (3.17) requires the condition $(mn, q) = 1$. By the theory of Hecke operators (as in the case of (2.9)) the sum (3.17) can be arranged into a sum over primitive cusp forms of level dividing $q$ with coefficients $\omega_j^*$ satisfying

(3.21) $$\omega_j^* \gg (q|s_j|)^{-1-\varepsilon}.$$

In the case of continuous spectrum we know the Fourier coefficients $\varphi_{\mathfrak{a}}(n,s)$ quite explicitly. We compute them by using the Eisenstein series for the modular group
$$E(z,s) = \tfrac{1}{2} y^s \sum\sum_{(c,d)=1} |cz+d|^{-2s}.$$

This has the Fourier expansion (3.3) with

(3.22) $$\varphi(n,s) = \pi^s \Gamma(s)^{-1} \zeta(2s) |n|^{-\frac{1}{2}} \eta(n,s)$$

where

(3.23) $$\eta(n,s) = \sum_{ad=|n|} (a/d)^{s-\frac{1}{2}}.$$

Every Eisenstein series $E_{\mathfrak{a}}(z,s)$ for the group $\Gamma = \Gamma_0(q)$ can be expressed as a linear combination of $E(dz,s)$ with $d|q$. Below we derive these representations.

Recall that $q$ is squarefree, so every cusp of $\Gamma_0(q)$ is equivalent to $\mathfrak{a} = 1/v$ with $v|q$. The complementary divisor $w = q/v$ is the width of $\mathfrak{a}$. We find that (by the arguments in [DI, p. 240], or [He, p. 534])
$$\sigma_{\mathfrak{a}}^{-1} \Gamma = \left\{ \begin{pmatrix} a/\sqrt{w} & b/\sqrt{w} \\ c\sqrt{w} & d\sqrt{w} \end{pmatrix} : \begin{pmatrix} a & b \\ c & d \end{pmatrix} \in \mathrm{SL}_2(\mathbb{Z}),\ c \equiv -av(q) \right\}.$$

Hence the cosets $\Gamma_\infty \backslash \sigma_{\mathfrak{a}}^{-1} \Gamma$ are parametrized by pairs of numbers $\{c\sqrt{w}, d\sqrt{w}\}$ with $(c,dw) = 1$ and $v|c$. Therefore the Eisenstein series for the cusp $\mathfrak{a} = 1/v$ is given by
$$E_{\mathfrak{a}}(z,s) = \sum_{\tau \in \Gamma_\infty \backslash \sigma_{\mathfrak{a}}^{-1} \Gamma} (\mathrm{Im}\,\tau z)^s = \frac{1}{2} \left(\frac{y}{w}\right)^s \sum\sum_{\substack{(c,dw)=1 \\ v|c}} |cz+d|^{-2s}.$$

Removing the condition $(d,c) = 1$ by Möbius inversion we get
$$E_{\mathfrak{a}}(z,s) = \frac{1}{2} \left(\frac{y}{w}\right)^s \sum_{(\delta,w)=1} \mu(\delta) \delta^{-2s} \sum\sum_{\substack{(c,d)\neq(0,0)\\(c,w)=1}} \left|\frac{cvz}{(\delta,v)} + d\right|^{-2s}$$
$$= \zeta(2s)^{-1} \sum_{(\delta,w)=1} \mu(\delta) \delta^{-2s} \sum_{\gamma|w} \mu(\gamma) \left(\frac{(\delta,v)}{\gamma q}\right)^s E\left(\frac{\gamma v z}{(\delta,v)}, s\right)$$
$$= \zeta_q(2s) \sum_{\delta|v}\sum_{\gamma|w} \mu(\delta\gamma)(\delta\gamma q)^{-s} E\left(\gamma \frac{v}{\delta} z, s\right)$$

where $\zeta_q(s)$ is the local zeta-function

(3.24) $$\zeta_q(s) = \prod_{p|q} (1-p^{-s})^{-1}.$$

Putting $v = \beta\delta$ we arrive at

$$(3.25) \qquad E_{\mathfrak{a}}(z,s) = \zeta_q(2s)\mu(v)(qv)^{-s} \sum_{\beta|v} \sum_{\gamma|w} \mu(\beta\gamma)\beta^s \gamma^{-s} E(\beta\gamma z, s).$$

By (3.22)–(3.25) we deduce that for $n \geq 1$, $(n,q) = 1$

$$(3.26) \qquad \varphi_{\mathfrak{a}}(n,s) = \frac{\mu(v)}{\Gamma(s)}\left(\frac{\pi}{qv}\right)^s \frac{\zeta_q(2s)}{\zeta(2s)} \frac{\eta(n,s)}{\sqrt{n}}.$$

Hence, for any cusp, $\eta_{\mathfrak{a}}(n,s) = \eta(n,s)$ if $(n,q) = 1$, and

$$(3.27) \qquad \omega_{\mathfrak{a}}(r) = \frac{4\pi}{qv} \frac{|\zeta_q(1+2ir)|^2}{|\zeta(1+2ir)|^2}.$$

Note that $\omega_{\mathfrak{a}}(r) \geq 4\pi w |\nu(q)\zeta(1+2ir)|^{-2}$, where $w$ is the width of the cusp and $\nu(q)$ is the index of the group.

## 4. Hecke $L$-functions

From now on we assume that $q$ is squarefree, odd. Let $\chi = \chi_q$ be the real, primitive character of conductor $q$; i.e., $\chi$ is given by the Jacobi-Legendre symbol

$$(4.1) \qquad \chi(n) = \left(\frac{n}{q}\right).$$

To any primitive form $f$ of level $q'|q$ we associate the $L$-functions

$$(4.2) \qquad L_f(s) = \sum_1^\infty \lambda_f(n) n^{-s} = \prod_{p \nmid q'} \left(1 - \lambda_f(p)p^{-s} + p^{-2s}\right)^{-1} \prod_{p|q'} \left(1 - \lambda_f(p)p^{-s}\right)^{-1}$$

and

$$(4.3) \qquad L_f(s,\chi) = \sum_1^\infty \lambda_f(n)\chi(n) n^{-s} = \prod_{p \nmid q}(1 - \lambda_f(p)\chi(p)p^{-s} + p^{-2s})^{-1}.$$

The latter is the $L$-function of the twisted form

$$f_\chi(z) = \sum_1^\infty a_f(n)\chi(n)n^{\frac{k-1}{2}} e(nz) \in S_k(\Gamma_0(q^2)).$$

Moreover, the completed $L$-function

$$(4.4) \qquad \Lambda_f(s,\chi) = \left(\frac{q}{2\pi}\right)^s \Gamma(s + \tfrac{k-1}{2}) L_f(s,\chi)$$

is entire and it satisfies the functional equation

$$(4.5) \qquad \Lambda_f(s,\chi) = w_f(\chi)\Lambda_f(1-s,\chi)$$

with $w_f(\chi) = \chi(-1)i^k$ (see Razar [R], for example). Note that $w_f(\chi)$ does not depend on $f$. From now on we assume that

(4.6) $$\chi(-1) = i^k$$

so that (4.5) holds with $w_f(\chi) = 1$ (otherwise all the central values $L_f(\frac{1}{2}, \chi)$ vanish).

*Remarks.* By a theorem of Winnie Li [L] the Euler product (4.3) and the functional equation (4.5) guarantee that $f_\chi$ is primitive. Also, we can see more explicitly the dependence of $w_f(\chi)$ on $f$ and $\chi$ as follows. If $f$ is primitive of level $q'|q$ then $w_f(\chi) = \chi(-1)\mu(q')\lambda_f(q')w_f$ and $w_f = i^k\mu(q')\lambda_f(q')$, so that $w_f(\chi) = \chi(-1)i^k$. Clearly all the above properties of $L_f(\frac{1}{2}, \chi)$ (including the definition (4.3)) remain true for any cusp form $f$ from the Hecke basis $\mathcal{F}$ (because the character $\chi$ kills the coefficients with $n$ not prime to $q$).

Using the functional equation (4.5) we shall represent the central values $L_f(\frac{1}{2}, \chi)$ by its partial sum of length about $O(kq)$. To this end we choose a function $G(s)$ which is holomorphic in $|\text{Re } s| \leq A$ such that

(4.7) $$G(s) = G(-s),$$
$$\Gamma(\tfrac{k}{2})G(0) = 1,$$
$$\Gamma(s + \tfrac{k}{2})G(s) \ll (|s|+1)^{-2A}$$

for some $A \geq 1$. Consider the integral

$$I = \frac{1}{2\pi i} \int_{(1)} \Lambda_f(s + \tfrac{1}{2}, \chi) G(s) s^{-1} \, ds.$$

Moving the integration to the line Re $s = -1$ and applying (4.5) we derive

$$\Lambda_f(\tfrac{1}{2}, \chi) G(0) = 2I.$$

On the other hand, integrating term by term, we derive

$$I = \sum_{1}^{\infty} \lambda_f(n) \chi(n) \left(\frac{q}{2\pi n}\right)^{\frac{1}{2}} V\left(\frac{n}{q}\right)$$

where $V(y)$ is the inverse Mellin transform of $(2\pi)^{-s} \Gamma(s + \tfrac{k}{2}) G(s) s^{-1}$,

(4.8) $$V(y) = \frac{1}{2\pi i} \int_{(1)} \Gamma(s + \tfrac{k}{2}) G(s) (2\pi y)^{-s} s^{-1} \, ds.$$

Hence by the normalization condition (4.7) we get:

LEMMA 4.1. *For any Hecke form $f \in \mathcal{F}$ we have*

(4.9) $$L_f(\tfrac{1}{2}, \chi) = 2 \sum_{1}^{\infty} \lambda_f(n) \chi(n) n^{-\frac{1}{2}} V(n/q).$$

Observe that $V(y)$ satisfies the following bounds

$$V(y) = 1 + O(y^A), \tag{4.10}$$

$$V(y) \ll (1+y)^{-A}, \tag{4.11}$$

$$V^{(\ell)}(y) \ll y^A(1+y)^{-2A}, \tag{4.12}$$

for $0 < \ell < A$ where the implied constant depends on that in (4.7). Actually $V(y)$ depends on the weight $k$. One can choose $G(s)$ depending on $k$ so that

$$V(y) \ll k(1+y/k)^{-A};$$

therefore the series (4.9) dies rapidly as soon as $n$ exceeds $kq$. If one is not concerned with the dependence of implied constants on the parameter $k$ then one has a simple choice $G(s) = \Gamma(k/2)^{-1}$ getting the incomplete gamma function

$$V(y) = \frac{1}{\Gamma(\frac{k}{2})} \int_{2\pi y}^{\infty} e^{-x} x^{\frac{k}{2}-1}\, dx.$$

By (2.5) and (4.9) we deduce that

$$L_f^2(\tfrac{1}{2},\chi) = 4 \sum_{(d,q)=1} d^{-1} \sum_{n_1} \sum_{n_2} \lambda_f(n_1 n_2) \frac{\chi(n_1 n_2)}{\sqrt{n_1 n_2}} V\left(\frac{dn_1}{q}\right) V\left(\frac{dn_2}{q}\right). \tag{4.13}$$

Now we have everything ready to begin working with the cubic moment of the central values $L_f(\tfrac{1}{2},\chi)$. From an analytic point of view it is natural to introduce the spectrally normalized cubic moment

$$\mathcal{C}_k(q) = \sum_{f \in \mathcal{F}} \omega_f L_f^3(\tfrac{1}{2},\chi) = \sum_{f \in \mathcal{F}^*} \omega_f^* L_f^3(\tfrac{1}{2},\chi) \tag{4.14}$$

where $\mathcal{F}$ is the Hecke orthonormal basis of $S_k(\Gamma_0(q))$. This differs from the arithmetically normalized cubic moment (1.2) by the coefficients $\omega_f^*$. Recall that the $\omega_f^*$ satisfy the lower bound (2.12). Therefore

$$\sum_{f \in \mathcal{F}^*} L_f^3(\tfrac{1}{2},\chi) \ll \mathcal{C}_k(q)(kq)^{1+\varepsilon} \tag{4.15}$$

for any $\varepsilon > 0$, where the implied constant depends only on $\varepsilon$. Hence for Theorem 1.1 we need to show that

$$\mathcal{C}_k(q) \ll q^\varepsilon \tag{4.16}$$

where the implied constant depends on $\varepsilon$ and $k$.

Applying (4.9) and (4.13) we write (4.14) as follows:

$$\mathcal{C}_k(q) = 8 \sum_{f \in \mathcal{F}} \omega_f \sum_n \sum_{n_1} \sum_{n_2} \lambda_f(n) \lambda_f(n_1 n_2) \frac{\chi(n n_1 n_2)}{\sqrt{n n_1 n_2}} V\left(\frac{n}{q}, \frac{n_2}{q}, \frac{n_2}{q}\right) \tag{4.17}$$

where

$$V(x, x_1, x_2) = V(x) \sum_{(d,q)=1} d^{-1} V(dx_1) V(dx_2). \tag{4.18}$$

Next by the Petersson formula (2.9) this is transformed into

$$\mathcal{C}_k(q) = \mathcal{D} + \sum_{c \equiv 0 \pmod{q}} c^{-2} \mathcal{S}(c) \tag{4.19}$$

where $\mathcal{D}$ is the contribution of the diagonal terms given by

$$\mathcal{D} = 8 \sum\sum_{(dn_1n_2,q)=1} (dn_1n_2)^{-1} V\left(\frac{n_1 n_2}{q}\right) V\left(\frac{dn_1}{q}\right) V\left(\frac{dn_2}{q}\right), \tag{4.20}$$

and $\mathcal{S}(c)$ is the contribution of the Kloosterman sums of modulus $c$ given by

$$\mathcal{S}(c) = 8 \sum_n \sum_{n_1} \sum_{n_2} \chi(nn_1n_2) S(n, n_1n_2; c) J\left(\frac{2\sqrt{nn_1n_2}}{c}\right) V\left(\frac{n}{q}, \frac{n_1}{q}, \frac{n_2}{q}\right). \tag{4.21}$$

## 5. Maass $L$-functions

To any even cusp form $u_j$ in the Hecke basis $\mathcal{U}$ of $\mathcal{L}(\Gamma \backslash \mathbb{H})$ we associate the $L$-function

$$L_j(s, \chi) = \sum_1^\infty \lambda_j(n) \chi(n) n^{-s}. \tag{5.1}$$

This has the Euler product of the type (4.3). Moreover the completed $L$-function

$$\Lambda_j(s, \chi) = \left(\frac{q}{\pi}\right)^s \Gamma\left(\frac{s + it_j}{2}\right) \Gamma\left(\frac{s - it_j}{2}\right) L_j(s, \chi), \tag{5.2}$$

is entire and it satisfies the functional equation

$$\Lambda_j(s, \chi) = \Lambda_j(1 - s, \chi). \tag{5.3}$$

Hence arguing as in Lemma 4.1 we deduce:

LEMMA 5.1. *For any even cusp form $u_j \in \mathcal{U}$,*

$$L_j(\tfrac{1}{2}, \chi) = 2 \sum_1^\infty \lambda_j(n) \chi(n) n^{-\frac{1}{2}} V_j(n/q) \tag{5.4}$$

*with $V_j(y)$ given by*

$$V_j(y) = \frac{1}{2\pi i} \int_{(1)} \Gamma\left(\frac{s + it_j}{2}\right) \Gamma\left(\frac{s - it_j}{2}\right) G_j(s) (\pi y)^{-s} s^{-1} ds. \tag{5.5}$$

*where $G_j(s)$ is any holomorphic function in $|\operatorname{Re} s| \leqslant A$ such that*

$$(5.6) \qquad G_j(s) = G_j(-s),$$

$$(5.7) \qquad \Gamma\left(\frac{1}{4} + \frac{it_j}{2}\right)\Gamma\left(\frac{1}{4} - \frac{it_j}{2}\right)G_j(0) = 1,$$

$$(5.8) \qquad \Gamma\left(\frac{s + it_j}{2}\right)\Gamma\left(\frac{s - it_j}{2}\right)G_j(s) \ll (|s| + 1)^{-2A}.$$

Observe that $V_j(y)$ satisfies the bounds (4.10)–(4.12).

To the Eisenstein series $E_\mathfrak{a}(z, \frac{1}{2} + ir)$ we associate the $L$-function

$$(5.9) \qquad L_{\mathfrak{a},r}(s,\chi) = \sum_1^\infty \eta_\mathfrak{a}(n, \tfrac{1}{2} + ir)\chi(n)n^{-s}.$$

It turns out that the $L_{\mathfrak{a},r}(s,\chi)$ is the same one for every cusp, indeed it is the product of two Dirichlet $L$-functions (see (3.26))

$$(5.10) \qquad L_{\mathfrak{a},r}(s,\chi) = L(s + ir, \chi)L(s - ir, \chi).$$

This satisfies the functional equation (5.3) (which can also be verified directly using the functional equation for $L(s,\chi)$; see [Da]), so (5.4) becomes

$$(5.11) \qquad |L(\tfrac{1}{2} + ir, \chi)|^2 = 2\sum_1^\infty \left(\sum_{ad=n}(a/d)^{ir}\right)\chi(n)n^{-\frac{1}{2}}V_r(n/q)$$

where $V_r(y)$ is given by the integral (5.5) with $t_j$ replaced by $r$ in (5.5)–(5.8).

Now we are ready to introduce the spectrally normalized cubic moment of the central values of $L$-functions associated with the even cusp forms and the Eisenstein series

$$(5.12) \qquad \mathcal{C}'_h(q) = {\sum_j}' h(t_j)\omega_j L_j^3(\tfrac{1}{2}, \chi) + \frac{1}{4\pi}\int_{-\infty}^\infty h(r)\omega(r)|L(\tfrac{1}{2} + ir)|^6\,dr$$

where the coefficients $\omega_j$ are given by (3.18) and $\omega(r) = \sum_\mathfrak{a}\omega_\mathfrak{a}(r)$. By (3.27) we obtain

$$(5.13) \qquad \omega(r) = \frac{4\pi}{q}\prod_{p|q}(1 + \frac{1}{p})\frac{|\zeta_q(1 + 2ir)|^2}{|\zeta(1 + 2ir)|^2}.$$

Note that the largest contribution to the continuous spectrum comes from the cusp of the largest width (which is the cusp zero). Now,

$$(5.14) \qquad \omega(r) \gg r^2((r^2 + 1)q)^{-1-\varepsilon}.$$

Assuming $h(r) \geq 0$ and $h(r) \geq 1$ if $-R \leq r \leq R$, we derive by (3.21) and (5.14) that

$$(5.15) \qquad {\sum_{|t_j|\leq R}}^\star L_j^3(\tfrac{1}{2}, \chi) + \int_{-R}^R |L(\tfrac{1}{2} + ir, \chi)|^6 \ell(r)\,dr \ll \mathcal{C}'_h(q)(Rq)^{1+\varepsilon}$$

for any $\varepsilon > 0$, where the implied constant depends only on $\varepsilon$. Therefore, for Theorem 1.4 we need to show that

$$\mathcal{C}'_h(q) \ll R^A q^\varepsilon \tag{5.16}$$

where the implied constant depends on $\varepsilon$ and the test function $h$.

Applying (5.4) and (5.11) we transform $\mathcal{C}'_h(q)$ by the Kuznetsov formula (3.17) into

$$\mathcal{C}'_h(q) = \mathcal{D}' + \sum_{c \equiv 0 (\mathrm{mod}\, q)} c^{-2} \mathcal{S}'(c), \tag{5.17}$$

where $\mathcal{D}'$ is the contribution of the diagonal terms and $\mathcal{S}'(c)$ is the contribution of the Kloosterman sums to modulus $c$. Here $\mathcal{D}'$ is similar to $\mathcal{D}$ in (4.20) and $\mathcal{S}'(c)$ is similar to $\mathcal{S}(c)$ in (4.21). From this point on our treatments of $\mathcal{C}_k(q)$ and $\mathcal{C}'_h(q)$ are almost identical except for technical details. Both cases are based on the same properties of the involved functions $J(x), J^+(x)$ and $J^-(x)$. Therefore, for notation economy, we choose to proceed further only with the classical cusp forms, i.e. we shall complete the proof of (4.16) and claim (5.16) by parallel arguments. The reader should note that the dependency on the spectral parameter in our estimates is polynomial at each step; hence the factor $R^A$ in (5.16) and (1.10).

## 6. Evaluation of the diagonal terms

Recall that $\mathcal{D}$ is the contribution to the cubic moment $\mathcal{C}_k(q)$ of the diagonal terms which is given by (4.20). Using the bounds (4.10)–(4.12) one shows that

$$\mathcal{D} \sim 8\Big(\frac{\varphi(q)}{q}\Big)^3 \int_1^\infty \int_1^\infty \int_1^\infty V\Big(\frac{yy_1}{q}\Big) V\Big(\frac{yy_2}{q}\Big) V\Big(\frac{y_1 y_2}{q}\Big) (yy_1 y_2)^{-1} \, dy \, dy_1 \, dy_2$$

$$\sim 16\Big(\frac{\varphi(q)}{q}\Big)^3 \int_1^q \Big(\int_1^{y_1} \Big(\int_1^{q/y_1} \frac{dy}{y}\Big) \frac{dy_2}{y_2}\Big) \frac{dy_1}{y_1} = \frac{8}{3}\Big(\frac{\varphi(q)}{q} \log q\Big)^3.$$

A more precise asymptotic can be derived from the complex integral

$$\mathcal{D} = 8(2\pi i)^{-3} \int_{(\varepsilon)} \int_{(\varepsilon)} \int_{(\varepsilon)} \Big(\frac{q}{2\pi}\Big)^{s+s_1+s_2}$$

$$\cdot \Gamma(s + \tfrac{k}{2})\Gamma(s_1 + \tfrac{k}{2})\Gamma(s_2 + \tfrac{k}{2}) G(s)G(s_1)G(s_2)$$

$$\cdot \zeta_q(s + s_1 + 1)\zeta_q(s + s_2 + 1)\zeta_q(s_1 + s_2 + 1)(s s_1 s_2)^{-1} \, ds \, ds_1 \, ds_2;$$

however, we only need the upper bound

$$\mathcal{D} \ll q^\varepsilon. \tag{6.1}$$

## 7. A partition of sums of Kloosterman sums

Recall that $\mathcal{S}(c)$ is the sum of Kloosterman sums $S(n, n_1 n_2; c)$ to the modulus $c \equiv 0 (\bmod q)$ given by (4.21) where $n, n_1, n_2$ run over all positive integers. The special function $J(x) = 4\pi i^k x^{-1} J_{k-1}(2\pi x)$ for $x = 2\sqrt{nn_1 n_2}/c$ which is attached to the Kloosterman sum $S(n, n_1 n_2; c)$ in (4.21) is itself a continuous analog of the latter. But this analogy is merely visual. One can compute asymptotically $J(2\sqrt{nn_1 n_2}/c)$ by the stationary phase method while the Kloosterman sum $S(n, n_1 n_2; c)$ requires more advanced arguments from algebraic geometry (see §§13 and 14).

All we need to know about the $J$-function is that it can be written as (see [W, p. 206])

$$J(x) = \operatorname{Re} W(x) e(x) \tag{7.1}$$

where $W(x)$ is a smooth function whose derivatives satisfy the bound (assuming $k \geqslant 12$)

$$x^\ell W^{(\ell)}(x) \ll x^{10}(1+x^2)^{-23/4} \tag{7.2}$$

for all $\ell \geqslant 0$ where the implied constant depends on $k$ and $\ell$. One could display the dependence on $k$, but we abandon this feature for the sake of notational simplicity.

To get hold on the variables $n, n_1, n_2$ of summation in (4.17) we split the range by a smooth partition of unity whose constituents are supported in dyadic boxes

$$\mathcal{N} = \left\{ (x, x_1, x_2) : 1 \leqslant \frac{x}{N}, \frac{x_1}{N_1}, \frac{x_2}{N_2} \leqslant 2 \right\} \tag{7.3}$$

with $N, N_1, N_2 \geqslant \frac{1}{2}$. Accordingly $\mathcal{S}(c)$ splits into sums of type

$$\mathcal{S}(W; c) = \sum_n \sum_{n_1} \sum_{n_2} \chi(nn_1 n_2) S(n, n_1 n_2; c) e(2\sqrt{nn_1 n_2}/c) W(n, n_1, n_2; c) \tag{7.4}$$

where

$$W(x, x_1, x_2; c) = P(x, x_1, x_2) V\left(\frac{x}{q}, \frac{x_1}{q}, \frac{x_2}{q}\right) W(2\sqrt{xx_1 x_2}/c) \tag{7.5}$$

and the $P$-functions are the constituents of the partition of unity. Precisely,

$$\mathcal{S}(c) = 8\operatorname{Re} \sum_{\mathcal{N}} \mathcal{S}(W; c). \tag{7.6}$$

We shall treat the sums of Kloosterman sums (7.4) in full generality. All we need to know about the function $W(x, x_1, x_2; c)$ is that it is smooth, supported in the box (7.3), and its partial derivatives satisfy the bound

$$(7.7) \qquad x^\ell x_1^{\ell_1} x_2^{\ell_2} |W^{(\ell,\ell_1,\ell_2)}(x, x_1, x_2; c)| \leqslant Q$$

for $0 \leqslant \ell, \ell_1, \ell_2 \leqslant A$ ($A$ is a large constant) with some $Q > 0$.

We also restrict the modulus $c$ to a dyadic segment

$$(7.8) \qquad C \leqslant c \leqslant 2C, \qquad c \equiv 0 \pmod{q},$$

with $C \geqslant q$. Notice that our particular function (7.5) satisfies (7.7) with

$$(7.9) \qquad Q \ll \left(\frac{NN_1N_2}{C^2}\right)^5 \left(1 + \frac{NN_1N_2}{C^2}\right)^{-23/4}$$

$$\cdot \left(1 + \frac{N}{q}\right)^{-A} \left(1 + \frac{N_1}{q}\right)^{-A} \left(1 + \frac{N_2}{q}\right)^{-A} \log q,$$

by virtue of (4.12) and (7.2) (the factor $\log q$ appears from the summation in $d$ in (4.18)).

*Remarks.* Using individual estimates for the Kloosterman sums in (7.4) one obtains

$$(7.10) \qquad \mathcal{S}(W; c) \ll cQNN_1N_2.$$

This bound is satisfactory for large $c$, but it is not sufficient in all ranges. In order to improve (7.10) one has to exploit some cancellation of the terms in (7.4), which is due to the variation in the argument of the twisted Kloosterman sum

$$\chi(nn_1n_2) S(n, n_1n_2; c) e(2\sqrt{nn_1n_2}/c)$$

with respect to $n, n_1, n_2$ (we shall get an extra cancellation by summing over $c$ as well).

## 8. Completing the sum $\mathcal{S}(W; c)$

The character $\chi(nn_1n_2)$ and the Kloosterman sum $S(n, n_1n_2; c)$ with $c \equiv 0 \pmod{q}$ in (7.4) are periodic in $n, n_1, n_2$ of period $c$ (however the exponential factor $e(2\sqrt{nn_1n_2}/c)$ is not). Thus splitting into residue classes and applying the Poisson summation formula for each class we obtain

$$(8.1) \qquad \mathcal{S}(W; c) = \sum_m \sum_{m_1} \sum_{m_2} G(m, m_1, m_2; c) \check{W}(m, m_1, m_2; c)$$

where $m, m_1, m_2$ run over all integers (we call the variables $m, m_1, m_2$ the "dual" of the "original" $n, n_1, n_2$), $G$ is the complete sum of Kloosterman sums twisted by characters to modulus $c$,

$$G(m, m_1, m_2; c) = \sum\sum\sum_{a,a_1,a_2 (\text{mod } c)} \chi(aa_1a_2)S(a, a_1a_2; c)e_c(am + a_1m_1 + a_2m_2) \tag{8.2}$$

(recall the standard notation $e_c(x) = e^{2\pi i x/c}$), and $\check{W}$ is given by the following Fourier integral:

$$\check{W}(m, m_1, m_2; c) \tag{8.3}$$
$$= \int_{\mathbb{R}^3} W(cx, cx_1, cx_2; c)e(2\sqrt{cxx_1x_2} - mx - m_1x_1 - m_2x_2) \, dx dx_1 dx_2.$$

Before employing advanced arguments we cut the series (8.1) using crude estimates for the sum (8.2) and the integral (8.3). By (10.16) we have

$$|G(m, m_1, m_2; c)| \leqslant c^3. \tag{8.4}$$

By (8.3) one gets directly $\check{W} \ll QNN_1N_2C^{-3}$; however on integrating by parts we can improve this bound to

$$\check{W}(m, m_1, m_2; c) \tag{8.5}$$
$$\ll QNN_1N_2C^{-3}\left(1 + \frac{|m|}{M}\right)^{-A}\left(1 + \frac{|m_1|}{M_1}\right)^{-A}\left(1 + \frac{|m_2|}{M_2}\right)^{-A}$$

where $MN = M_1N_1 = M_2N_2 = C + \sqrt{NN_1N_2} = D$, say. Since $A$ is a large constant it shows that $\check{W}$ is very small outside the box

$$|m| \leqslant Mq^\varepsilon, \ |m_1| \leqslant M_1q^\varepsilon, \ |m_2| \leqslant M_2q^\varepsilon. \tag{8.6}$$

We denote this box by $\mathcal{M}$ (and say $\mathcal{M}$ is "dual" to $\mathcal{N}$). Estimating the tail of the series (8.1) with $(m, m_1, m_2) \notin \mathcal{M}$ by using (8.4) and (8.5) we are left with

$$\mathcal{S}(W; c) = \sum\sum\sum_{(m,m_1,m_2) \in \mathcal{M}} G(m, m_1, m_2; c)\check{W}(m, m_1, m_2; c) + O(QNN_1N_2q^{-3}). \tag{8.7}$$

Our next goal is to pull out from $\check{W}(m, m_1, m_2; c)$ the phase factor $e(-mm_1m_2/c)$ and then to separate the variables in

$$K(m, m_1, m_2; c) = e\left(\frac{mm_1m_2}{c}\right)\check{W}(m, m_1, m_2; c). \tag{8.8}$$

We choose the method of Fourier transform (because it is in harmony with the forthcoming computations of the sum $G(m, m_1, m_2; c)$ in §10); however the Mellin transform would do the job as well.

Throughout $D^{\mathbf{a}}$ stands for the differential operator

$$D^{\mathbf{a}}F(x, x_1, x_2) = x^a x_1^{a_1} x_2^{a_2} F^{(a,a_1,a_2)}(x, x_1, x_2)$$

where $F^{(a,a_1,a_2)}$ is the partial derivative of order $\mathbf{a} = (a, a_1, a_2)$.

LEMMA 8.1. *Let $W(x, x_1, x_2)$ be a smooth function supported in the dyadic box*

(8.9) $$\mathcal{X} = \left\{(x, x_1, x_2) : 1 \leqslant \frac{x}{X}, \frac{x_1}{X_1}, \frac{x_2}{X_2} \leqslant 2\right\}$$

*with $X, X_1, X_2 > 0$, and its partial derivatives satisfy the bound*

(8.10) $$|D^{\mathbf{a}}W(x, x_1, x_2)| \leqslant 1$$

*for any $\mathbf{a} = (a, a_1, a_2) \leqslant 18$ component by component. Put*

(8.11) $$L(x, x_1, x_2) = W(x, x_1, x_2)e(2\sqrt{xx_1x_2}),$$
(8.12) $$K(y, y_1, y_2) = \hat{L}(y, y_1, y_2)e(yy_1y_2),$$

*where $\hat{L}$ is the Fourier transform of $L$. Then the $L_1$-norm of the Fourier transform of $K$ satisfies*

(8.13) $$\int_{\mathbb{R}^3} |\hat{K}(t, t_1, t_2)| dt dt_1 dt_2 \ll (XX_1X_2)^{\frac{1}{4}} \log(1 + XX_1X_2) + (XX_1X_2)^{-2}$$

*where the implied constant is absolute.*

*Proof.* Our arguments are direct but quite long as we go through the Fourier inversion several times to avoid problems with the stationary phase. Define $Y, Y_1, Y_2$ and $Z$ by

(8.14) $$XX_1X_2 = Z^2, \ XY = X_1Y_1 = X_2Y_2 = Z.$$

Note that $YY_1Y_2 = Z$. First we estimate the derivatives of

(8.15) $$\hat{L}(y, y_1, y_2) = \int_{\mathbb{R}^3} W(x, x_1, x_2)e(2\sqrt{xx_1x_2} - xy - x_1y_1 - x_2y_2)d\mathbf{x}$$

at any point $\mathbf{y} = (y, y_1, y_2) \in \mathbb{R}^3$ outside the box

(8.16) $$\mathcal{Y} = \left\{(y, y_1, y_2) : \frac{1}{3} \leqslant \frac{y}{Y}, \frac{y_1}{Y_1}, \frac{y_2}{Y_2} \leqslant 3\right\}.$$

One may say that $\mathcal{Y}$ is "dual" to $\mathcal{X}$. Applying the operator $D^{\mathbf{a}}$ with respect to the variables $(y, y_1, y_2)$ we get

$$D^{\mathbf{a}}\hat{L}(y, y_1, y_2) = \int_{\mathbb{R}^3} W(x, x_1, x_2)(-2\pi ixy)^a(-2\pi ix_1y_1)^{a_1}(-2\pi ix_2y_2)^{a_2}$$
$$\cdot e(2\sqrt{xx_1x_2} - xy - x_1y_1 - x_2y_2)d\mathbf{x}.$$

Hence $D^{\mathbf{a}}\hat{L}(y, y_1, y_2) \ll Z^2(|y|X)^a(|y_1|X_1)^{a_1}(|y_2|X_2)^{a_2}$ by trivial estimation; however we can do better by partial integration. Since $\mathbf{y} = (y, y_1, y_2) \notin \mathcal{Y}$ we may assume without loss of generality that $y$ is not in the segment $\frac{1}{3}Y \leqslant y \leqslant 3Y$, so $xy$ does not match $2\sqrt{xx_1x_2}$. If $Z + |y|X > 1$ then on integrating by parts with respect to $x$ eighteen times we gain the factor $(Z + |y|X)^{18}$;

otherwise we do not integrate by parts with respect to $x$ at all. In any case we save the factor $(1 + Z + |y|X)^{18}$. The same operation can be applied simultaneously with respect to the other variables $x_1, x_2$, provided $y_1, y_2$ do not satisfy $\frac{1}{3}Y_i \leqslant y_i \leqslant 3Y_i$, gaining the factors $(1+Z+|y_i|X_i)^{18}$. If one or both variables $y_1, y_2$ do satisfy $\frac{1}{3}Y_i \leqslant y_i \leqslant 3Y_i$ then we do not integrate by parts with respect to $x_i$, but still claim the factor $(1 + Z + |y_i|X_i)^6$ by borrowing it from the gain in the variable $x$. Thus for any $\mathbf{y} = (y, y_1, y_2) \in \mathbb{R}^3$ outside the box $\mathcal{Y}$ and any $\mathbf{a} = (a, a_1, a_2)$ we have

(8.17)
$$D^{\mathbf{a}}\hat{L}(y, y_1, y_2) \ll Z^2(|y|X)^a(|y|_1|X_1)^{a_1}(|y_2|X_2)^{a_2}$$
$$\cdot (1 + Z + |y|X)^{-6}(1 + Z + |y_1|X_1)^{-6}(1 + Z + |y_2|X_2)^{-6}.$$

Applying $D^{\mathbf{a}}$ to (8.12) we derive by (8.17) that

$$D^{\mathbf{a}}K(y, y_1, y_2) \ll Z^2(|y|X + |yy_1y_2|)^a(|y_1|X_1 + |yy_1y_2|)^{a_1}(|y_2|X_2 + |yy_1y_2|)^{a_2}$$
$$\cdot (1 + Z + |y|X)^{-6}(1 + Z + |y_1|X_1)^{-6}(1 + Z + |y_2|X_2)^{-6}.$$

Here we have $|y|X + |yy_1y_2| \leqslant |y|XZ^{-2}(1 + Z + |y_1|X_1)(1 + Z + |y_2|X_2)$ and similar inequalities hold for the other two combinations. Therefore if $\mathbf{a} \leqslant 2$ and $\mathbf{y} \notin \mathcal{Y}$,

(8.18)  $$D^{\mathbf{a}}K(y, y_1, y_2) \ll Z^2(|y|XZ^{-2})^a(|y_1|X_1Z^{-2})^{a_1}(|y_2|X_2Z^{-2})^{a_2}$$
$$\cdot (1 + |y|X)^{-2}(1 + |y_1|X_1)^{-2}(1 + |y_2|X_2)^{-2}.$$

Now we proceed to the estimation of $D^{\mathbf{a}}K(\mathbf{y})$ in the box $\mathcal{Y}$. In this range before integrating by parts we pull out the exponential factor $e(-yy_1y_2)$ from the Fourier integral (8.15). To this end we arrange the amplitude function $2\sqrt{xx_1x_2} - xy - x_1y_1 - x_2y_2$ in the following form

$$-yy_1y_2 + \frac{1}{y}(x_1 - yy_2)(x_2 - yy_1) - \frac{1}{y}(y\sqrt{x} - \sqrt{x_1x_2})^2.$$

Introducing this into (8.15) and changing the variables of integration $\mathbf{x} = (x, x_1, x_2)$ into $\mathbf{v} = (v, v_1, v_2)$ by the formulas $x_1 = (v_1 + y_2\sqrt{y})\sqrt{y}$, $x_2 = (v_2 + y_1\sqrt{y})\sqrt{y}$ and $x = (v + \sqrt{x_1x_2/y})^2/y$ we get

(8.19)  $$K(\mathbf{y}) = \int_{\mathbb{R}^3} H(\mathbf{v}; \mathbf{y})e(-v^2 + v_1v_2)\,dv\,dv_1\,dv_2$$

where the kernel function is given by

(8.20)  $$H(\mathbf{v}; \mathbf{y}) = 2(x(\mathbf{v}; \mathbf{y})y)^{\frac{1}{2}}W(x(\mathbf{v}; \mathbf{y}), (v_1 + y_2\sqrt{y})\sqrt{y}, (v_2 + y_1\sqrt{y})\sqrt{y})$$

with

$$(8.21) \qquad x(\mathbf{v};\mathbf{y}) = [v + (v_1 + y_2\sqrt{y})^{\frac{1}{2}}(v_2 + y_1\sqrt{y})^{\frac{1}{2}}]^2 y^{-1}.$$

Notice that the new variables $\mathbf{v} = (v, v_1, v_2)$ are restricted by

$$(8.22) \qquad \begin{aligned} v\sqrt{y} &= y\sqrt{x} - \sqrt{x_1 x_2} \ll YX^{\frac{1}{2}} + (X_1 X_2)^{\frac{1}{2}} = 2YX^{\frac{1}{2}}, \\ v_1\sqrt{y} &= x_1 - yy_2 \ll X_1 + YY_2 = 2X_1, \\ v_2\sqrt{y} &= x_2 - yy_2 \ll X_2 + YY_1 = 2X_2, \end{aligned}$$

by virtue of the support of $W(x, x_1, x_2)$ being in the box $\mathcal{X}$ and the point $(y, y_1, y_2)$ being in the box $\mathcal{Y}$. Put

$$(8.23) \qquad p(\mathbf{v}) = -\int_0^v \int_0^{v_1} \int_0^{v_2} e(-u^2 + u_1 u_2) du\, du_1\, du_2.$$

Clearly, we have $p(\mathbf{v}) \ll \log(1 + |v_1 v_2|)$ on $\mathbb{R}^3$. Integrating (8.19) by parts we get

$$K(\mathbf{y}) = \int_{\mathbb{R}^3} p(\mathbf{v}) \frac{\partial^3 H(\mathbf{v};\mathbf{y})}{\partial v\, \partial v_1\, \partial v_2} dv\, dv_1\, dv_2.$$

Applying the operator $D^{\mathbf{a}}$ and the restrictions (8.22) we derive the following estimate

$$(8.24) \qquad D^{\mathbf{a}} K(\mathbf{y}) \ll \left| \frac{\partial^3 D^{\mathbf{a}} H(\mathbf{v};\mathbf{y})}{\partial v\, \partial v_1\, \partial v_2} \right| Z^{\frac{3}{2}} \log(1 + Z)$$

for some $\mathbf{v} = (v, v_1, v_2)$. We need (8.24) for all $\mathbf{y} \in \mathcal{Y}$ and $\mathbf{a} \leqslant 2$; therefore we have to differentiate $H(\mathbf{v};\mathbf{y})$ up to nine times. However we only show details for the first order partial derivatives, the higher order ones, being estimated by repeating the arguments, are left for checking to careful readers.

We begin by the following estimate

$$(8.25) \qquad H(\mathbf{v};\mathbf{y}) \ll (XY)^{\frac{1}{2}} = Z^{\frac{1}{2}}.$$

Next we estimate the partial derivatives of $x(\mathbf{v};\mathbf{y})$ with respect to $v, v_1, v_2$ and $y, y_1, y_2$ in the range restricted by the support of $W(x, x_1, x_2)$. By (8.21) we derive the following estimates

$$\frac{\partial x}{\partial v} = 2[v + (v_1 + y_2\sqrt{y})^{\frac{1}{2}}(v_2 + y_1\sqrt{y})^{\frac{1}{2}}]y^{-1} \ll \left(\frac{X}{Y}\right)^{\frac{1}{2}} = XZ^{-\frac{1}{2}},$$

$$\frac{\partial x}{\partial v_1} = \left(\frac{v_2 + y_1\sqrt{y}}{v_1 + y_2\sqrt{y}}\right)^{\frac{1}{2}} \frac{[*]}{y} \ll \left(\frac{XX_2}{YX_1}\right)^{\frac{1}{2}} = X\left(\frac{X_2}{ZX_1}\right)^{\frac{1}{2}},$$

$$\frac{\partial x}{\partial v_2} = \left(\frac{v_1 + y_2\sqrt{y}}{v_2 + y_1\sqrt{y}}\right)^{\frac{1}{2}} \frac{[*]}{y} \ll \left(\frac{XX_1}{YX_2}\right)^{\frac{1}{2}} = X\left(\frac{X_1}{ZX_2}\right)^{\frac{1}{2}},$$

$$y\frac{\partial x}{\partial y} = \{(\ldots)^{-\frac{1}{2}}(\ldots)^{\frac{1}{2}}y_2 + (\ldots)^{\frac{1}{2}}(\ldots)^{-\frac{1}{2}}y_1\}[*]y^{-\frac{1}{2}} - [*]^2 y^{-1}$$

$$\ll \{(X_2/X_1)^{\frac{1}{2}}Y_2 + (X_1/X_2)^{\frac{1}{2}}Y_1\}X^{\frac{1}{2}} + X = 3X,$$

$$y_1\frac{\partial x}{\partial y_1} = y_1(\ldots)^{\frac{1}{2}}(\ldots)^{-\frac{1}{2}}[*]y^{-\frac{1}{2}} \ll Y_1(XX_1/X_2)^{\frac{1}{2}} = X,$$

$$y_2\frac{\partial x}{\partial y_2} = y_2(\ldots)^{-\frac{1}{2}}(\ldots)^{\frac{1}{2}}[*]y^{-\frac{1}{2}} \ll Y_2(XX_2/X_1)^{\frac{1}{2}} = X.$$

Continuing the differentiation along the above lines one shows that

(8.26) $$D^{\mathbf{a}}x(\mathbf{v};\mathbf{y}) \ll X$$

and

(8.27) $$\frac{\partial^{(\alpha,\alpha_1,\alpha_2)}}{\partial v^\alpha \partial v_1^{\alpha_1} \partial v_2^{\alpha_2}} D^{\mathbf{a}}x(\mathbf{v};\mathbf{y}) \ll X(X_2/X_1)^{\frac{\alpha_1}{2}}(X_1/X_2)^{\frac{\alpha_2}{2}} Z^{-\frac{1}{2}(\alpha+\alpha_1+\alpha_2)}$$

if $(\alpha, \alpha_1, \alpha_2) \leqslant 1$ and $\mathbf{a} = (a, a_1, a_2) \leqslant 2$. Now we are ready to estimate the partial derivatives of $H(\mathbf{v};\mathbf{y})$. By (8.20) we obtain the following estimates

$$\frac{\partial H}{\partial v} = \left(\frac{y}{x}\right)^{\frac{1}{2}}\left(W + 2x\frac{\partial W}{\partial x}\right)\frac{\partial x}{\partial v} \ll 1,$$

$$\frac{\partial H}{\partial v_1} = \left(\frac{y}{x}\right)^{\frac{1}{2}}\left(W + 2x\frac{\partial W}{\partial x}\right)\frac{\partial x}{\partial v_1} + 2y\sqrt{x}\frac{\partial W}{\partial x_1}$$

$$\ll (X_2/X_1)^{\frac{1}{2}} + YX^{\frac{1}{2}}X_1^{-1} = 2(X_2/X_1)^{\frac{1}{2}},$$

$$\frac{\partial H}{\partial v_2} \ll (X_1/X_2)^{\frac{1}{2}} \quad \text{(by interchanging variables)},$$

$$y\frac{\partial H}{\partial y} = y\left(\frac{\partial x}{\partial y}\right)(\frac{y}{x})^{\frac{1}{2}}\left(W + x\frac{\partial W}{\partial x}\right) \ll Z^{\frac{1}{2}},$$

$$y_1\frac{\partial H}{\partial y_1} = y_1\left(\frac{\partial x}{\partial y_1}\right)(\frac{y}{x})^{\frac{1}{2}}\left(W + (xy)^{\frac{1}{2}}\frac{\partial W}{\partial x}\right) + 2y_1 y(xy)^{\frac{1}{2}}\frac{\partial W}{\partial x_2} \ll Z^{\frac{1}{2}},$$

$$y_2\frac{\partial H}{\partial y_2} \ll Z^{\frac{1}{2}} \quad \text{(by interchanging variables)}.$$

Continuing the differentiation along the above lines one shows that

(8.28) $$D^{\mathbf{a}} H(\mathbf{v}; \mathbf{y}) \ll Z^{\frac{1}{2}},$$

(8.29) $$\frac{\partial^{(\alpha,\alpha_1,\alpha_2)}}{\partial v^\alpha \partial v_1^{\alpha_1} \partial v_2^{\alpha_2}} D^{\mathbf{a}} H(\mathbf{v}; \mathbf{y}) \ll \left(\frac{X_2}{X_1}\right)^{\frac{\alpha_1}{2}} \left(\frac{X_1}{X_2}\right)^{\frac{\alpha_2}{2}} Z^{\frac{1}{2}(1-\alpha-\alpha_1-\alpha_2)}$$

if $(\alpha, \alpha_1, \alpha_2) \leq 1$ and $\mathbf{a} = (a, a_1, a_2) \leq 2$. In particular,

$$\frac{\partial^3 D^{\mathbf{a}} H(\mathbf{v}; \mathbf{y})}{\partial v \partial v_1 \partial v_2} \ll Z^{-1}.$$

Inserting this into (8.24) we conclude that for $\mathbf{a} \leq 2$ and $\mathbf{y} \in \mathcal{Y}$

(8.30) $$D^{\mathbf{a}} K(\mathbf{y}) \ll Z^{\frac{1}{2}} \log(1+Z).$$

Finally we are ready to complete the proof of Lemma 8.1. Combining (8.18) and (8.30) we deduce that the Fourier transform of $K(y, y_1, y_2)$ satisfies

(8.31) $$\hat{K}(t, t_1, t_2) \ll (1+|t|X_1 X_2)^{-2}(1+|t_1|XX_2)^{-2}(1+|t_2|XX_1)^{-2}$$
$$+ Z^{\frac{3}{2}}(1+|t|Y)^{-2}(1+|t_1|Y_1)^{-2}(1+|t_2|Y_2)^{-2}$$

where the first term comes by partial integration outside the box $\mathcal{Y}$ by (8.18), and the second term comes by partial intergration in the box $\mathcal{Y}$ by (8.30) (needless to say the transition through the boundary of $\mathcal{Y}$ is made with a smooth partition of unity to avoid boundary terms in partial integration). By (8.31) we obtain

$$\int_{\mathbb{R}^3} |\hat{K}(t, t_1, t_2)| dt dt_1 dt_2 \ll Z^{-4} + Z^{\frac{1}{2}} \log(1+Z)$$

which is the bound (8.13).

*Remarks.* The term $Z^{-4}$ is significant only if $Z \leq 1$, it could be improved if we dealt with the factor $e(yy_1 y_2)$ in passing from (8.17) to (8.18) by stationary phase methods rather than by direct differentiation.

By Lemma 8.1 we write

(8.32) $$\check{W}(m, m_1, m_2; c) = e\left(-\frac{mm_1 m_2}{c}\right) K(m, m_1, m_2; c)$$

where for $c$ fixed $K(y, y_1, y_2; c) = K(y, y_1, y_2)$ is a smooth function on $\mathbb{R}^3$ whose Fourier transform satisfies

(8.33) $$\int_{\mathbb{R}^3} |\hat{K}(t, t_1, t_2)| dt dt_1 dt_2 \ll C^{-1}(Z^{\frac{1}{2}} \log(1+Z) + Z^{-4})$$

with $Z = \sqrt{NN_1 N_2}/C$. We shall see that the factor $e(-mm_1 m_2/c)$ which we extracted from $\check{W}(m, m_1, m_2; c)$ in (8.32) cancels out with the factor $e(mm_1m_2/c)$ which appears in the formula (10.1) for the character sum $G(m, m_1, m_2; c)$. Originally $G(m, m_1, m_2; c)$ was defined by (8.2); however, anticipating the cancellation of the above character, we introduce the modified sum

$$(8.34) \qquad G'(m, m_1, m_2; c) = e(-\frac{mm_1m_2}{c})G(m, m_1, m_2; c),$$

and we refer the reader to Lemma 10.2 to see another expression for $G'(m, m_1, m_2; c)$. Inserting (8.32), (8.34) into (8.7) and applying the Fourier inversion to $K(m, m_1, m_2; c)$ we deduce by (8.33) that

$$(8.35) \qquad \mathcal{S}(W; c) \ll NN_1N_2q^{-3} + C^{-1}(Z^{\frac{1}{2}}\log(1+Z) + Z^{-4})$$
$$\Big|\sum\sum\sum_{(m,m_1,m_2)\in\mathcal{M}} G'(m, m_1, m_2; c)\psi(m, m_1, m_2)\Big|$$

where $\psi(m, m_1, m_2) = e(-mt - m_1t_1 - m_2t_2)$ for some real numbers $t, t_1, t_2$.

It does not matter what the numbers $t, t_1, t_2$ are since in the following sections we are going to establish estimates for sums of type

$$(8.36)$$
$$\mathcal{G}(M, M_1, M_2; C) = \sum_{\substack{C \leqslant c < 2C \\ c \equiv 0 (\bmod q)}} c^{-2} \Big|\sum_m \sum_{m_1} \sum_{m_2} G'(m, m_1, m_2; c)\alpha_m \beta_{m_1, m_2}\Big|$$

for any complex coefficients $\alpha_m$ for $|m| \leqslant M$ and $\beta_{m_1, m_2}$ for $|m_1| \leqslant M_1$, $|m_2| \leqslant M_2$ which are bounded.

PROPOSITION 8.2. *Let $C \geqslant q > 1$ and $M, M_1, M_2 \geqslant 1$. Let $\alpha_m$ and $\beta_{m_1, m_2}$ be complex numbers with $|\alpha_m| \leqslant 1$ for $|m| \leqslant M$ and $|\beta_{m_1, m_2}| \leqslant 1$ for $|m_1| \leqslant M_1, |m_2| \leqslant M$. Then*

$$(8.37)$$
$$\mathcal{G}(M, M_1, M_2; C) \ll C\left(1 + (M + M_1 + M_2)q^{-1}\right)^2 \tau^2(q)$$
$$+ \left(1 + CMq^{-2}\right)^{\frac{1}{2}}\left(1 + M_1M_2q^{-1}\right)^{\frac{1}{2}}(CMM_1M_2)^{\frac{1}{2}+\varepsilon}$$

*for any $\varepsilon > 0$, the implied constant depending only on $\varepsilon$.*

## 9. Estimation of the cubic moment – Conclusion

Using Proposition 8.2 we are ready to finish the proof of Theorem 1.1. Recall that the cubic moment (4.14) was transformed by Petersson's formula to (4.19) where $\mathcal{D}$ is the contribution of the diagonal terms as estimated in (6.1). The other terms $\mathcal{S}(c)$ are sums of Kloosterman sums. These are partitioned into dyadic boxes in (7.6). We also divided the range of the modulus $c$ into dyadic segments (7.8). Now, given a box $\mathcal{N}$ of size $N \times N_1 \times N_2$ and $C \geqslant q$, we need to estimate

$$(9.1) \qquad \mathcal{T}(C) = \sum_{\substack{C \leqslant c < 2C \\ c \equiv 0 (\bmod q)}} c^{-2} \mathcal{S}(W; c).$$

To this end we apply (8.35) getting

$$(9.2) \qquad \mathcal{T}(C) \ll QNN_1N_2 C^{-1} q^{-4}$$
$$+ QC^{-1}(Z^{\frac{1}{2}} \log(1+Z) + Z^{-4}) \mathcal{G}(M, M_1, M_2; C)$$

where $Q$ satisfies (7.9) and $M \times M_1 \times M_2$ is the size of the "dual" box. Let us recall that $Z = \sqrt{NN_1N_2}/C$, and $MN = M_1N_1 = M_2N_2 = D$ with $D = C + \sqrt{NN_1N_2}$. Introducing (8.37) into (9.2) we get

$(9.3)$
$$\mathcal{T}(C) \ll QNN_1N_2 C^{-1} q^{-4} + Q(Z^{\frac{1}{2}} + Z^{-4})\left(1 + \frac{D}{qN} + \frac{D}{qN_1} + \frac{D}{qN_2}\right)^2 D^\varepsilon$$
$$+ QC^{-1}(Z^{\frac{1}{2}} + Z^{-4})\left(1 + \frac{CD}{q^2 N}\right)^{\frac{1}{2}} \left(1 + \frac{D^2}{qN_1N_2}\right)^{\frac{1}{2}} \left(\frac{CD^3}{NN_1N_2}\right)^{\frac{1}{2}} D^\varepsilon.$$

Taking into account the bound for $W$ given in (7.9) we see easily that the right-hand side of (9.3) multiplied by $C^\varepsilon$ is the largest for $C = \sqrt{NN_1N_2}$ and $N = N_1 = N_2 = q$, and it is bounded by $O(q^{3\varepsilon})$. Therefore

$$(9.4) \qquad \mathcal{T}(C) \ll q^{3\varepsilon} C^{-\varepsilon}$$

for any $C \geqslant q$ and any $\varepsilon > 0$, the implied constant depending on $\varepsilon$. Finally gathering together all the pieces (9.4) and (6.1) into (4.19) we get $\mathcal{C}_k(q) \ll q^{2\varepsilon}$ which completes the proof of Theorem 1.1.

Of course, we still have to prove Proposition 8.2 for which we spend the rest of this paper.

## 10. Evaluation of $G(m, m_1, m_2; c)$

Recall that $G(m, m_1, m_2; c)$ is defined by (8.2). Put $c = qr$. We begin by opening the Kloosterman sum

$$S(a, a_1 a_2; c) = \sum_{d (\bmod c)}^\star e\left(\frac{ad + a_1 a_2 \bar{d}}{c}\right).$$

The sum over $a(\bmod c)$ vanishes unless $d + m \equiv 0 (\bmod r)$ in which case it is equal to $r\tau(\chi)\chi((d+m)/r)$ where $\tau(\chi)$ is the Gauss sum

$$\tau(\chi) = \sum_{x(\bmod q)} \chi(x) e_q(x).$$

Similarly the sum over $a_1 \pmod c$ vanishes unless $a_2 \bar{d} + m_1 \equiv 0 \pmod r$ in which case it is equal to $r\tau(\chi)\chi(d)\chi((a_2 + dm_1)/r)$. Now the sum over $a_2$ is given by

$$\sum_{\substack{a_2 \pmod c \\ a_2 \equiv -dm_1 \pmod r}} \chi(a_2)\chi((a_2 + dm_1)/r) e(a_2 m_2/qr)$$

$$= \sum_{u \pmod q} \chi(ur - dm_1)\chi(u) e\left(\frac{um_2}{q} - \frac{dm_1 m_2}{qr}\right).$$

Finally we sum over $d \pmod c$ with $(d, c) = 1, d \equiv -m \pmod r$. These conditions imply that $m$ and $r$ are co-prime. Having recorded that $(m, r) = 1$ the condition $(d, c) = 1$ is redundant so we are free to run over $d \pmod c$ getting

$$\sum_{\substack{d \pmod c \\ d \equiv -m \pmod r}} \chi(d)\chi((d + m)/r)\chi(ur - dm_1) e(-dm_1 m_2/qr)$$

$$= \sum_{v \pmod q} \chi(v)\chi(vr - m)\chi(ur - (vr - m)m_1) e\left(\frac{-vm_1 m_2}{q} + \frac{mm_1 m_2}{qr}\right).$$

Gathering the above results we obtain (for $(m, r) = 1$)

(10.1) $\quad G(m, m_1, m_2; qr) = (r\tau(\chi))^2 e_{qr}(mm_1 m_2) H_r(m, m_1, m_2; q)$

where

$$H_r(m, m_1, m_2; q) = \sum\sum_{u,v \pmod q} \chi[uv(vr-m)(ur-(vr-m)m_1)] e_q(um_2 - vm_1 m_2).$$

Changing $u$ into $u + vm_1$ we arrive at

(10.2) $\quad H_r(m, m_1, m_2; q) = \sum\sum_{u,v \pmod q} \chi[v(u+vm_1)(vr-m)(ur+mm_1)] e_q(um_2).$

Clearly the character sum $H_r(m, m_1, m_2; q)$ is symmetric in $m_1, m_2$ by (10.1) and (8.2). Note also that it is multiplicative in $q$; precisely if $(q_1, q_2) = 1$ then

(10.3) $\quad H_r(m, m_1, m_2; q_1 q_2) = H_r(m, m_1 \bar{q}_1, m_2; q_2) H_r(m, m_1, m_2 \bar{q}_2; q_1)$

where $\bar{q}_1, \bar{q}_2$ are the multiplicative inverses of $q_1, q_2$ to moduli $q_2, q_1$ respectively. (The implied characters are $\chi_{q_1 q_2}, \chi_{q_2}$, and $\chi_{q_1}$ respectively.) This property reduces our problem of computing $H_r(m, m_1, m_2; q)$ to prime modulus.

LEMMA 10.1. *Suppose $q$ is prime. Then $H_r(m, m_1, m_2; q)$ is given by*

(10.4) $\quad\quad\quad\quad \chi^2(mm_1 m_2)\tau^2(\chi), \quad\text{if } q|r,$

(10.5) $$\frac{1}{q-1}R(m;q)R(m_1;q)R(m_2;q), \quad \text{if } q \nmid r \text{ and } q|mm_1m_2,$$

(10.6) $$H(\bar{r}mm_1m_2;q), \quad \text{if } q \nmid rmm_1m_2.$$

Here $\tau(\chi)$ is the Gauss sum, $R(m,q) = S(0,m;q)$ is the Ramanujan sum, and

(10.7) $$H(w;q) = \sum\sum_{u,v \pmod{q}} \chi(uv(u+1)(v+1))e_q((uv-1)w).$$

*Remark.* In the sequel we consider (10.7) as the definition of $H(w;q)$ for all $q$, not necessarily prime.

*Proof.* We use the expression (10.2). If $q|r$ then (10.2) reduces to

$$\chi(-m^2m_1)\sum_u\sum_v \chi(v(u+vm_1))e_q(um_2) = \chi^2(mm_1m_2)\tau^2(\chi).$$

Now suppose $q \nmid r$. Then changing $u,v$ into $u\bar{r}, v\bar{r}$ we get

(10.8) $$\sum_u\sum_v \chi[v(u+vm_1)(v-m)(u+mm_1)]e_q(u\bar{r}m_2).$$

If $q|m_1$ then (10.8) reduces to

$$\sum_u \chi^2(u)e_q(u\bar{r}m_2)\sum_v \chi(v(v-m)) = R(m_2;q)R(m;q)$$

which agrees with (10.5). The case $q|m_2$ follows by the symmetry. Now suppose $q \nmid rm_1m_2$. Then (10.8) becomes

(10.9) $$\sum_u\sum_v \chi[v(u+v)(v-m)(u+m)]e_q(u\bar{r}m_1m_2).$$

If $q|m$ then (10.9) reduces to

$$\sum_u\sum_v \chi[uv^2(u+v)]e_q(u\bar{r}m_1m_2) = 1$$

which agrees with (10.5). Finally, if $q \nmid rmm_1m_2$ then (10.9) becomes (10.6) by changing variables.

*Remarks.* Note that $H_r(m,m_1,m_2;q)$ vanishes if $q|(r,mm_1m_2)$, and it is equal to $\tau^2(\chi) = \chi(-1)q$ if $q|r, q \nmid mm_1m_2$. One can express the two cases (10.5) and (10.6) in one form

(10.10)
$$H_r(m,m_1,m_2;q) = H(\bar{r}mm_1m_2;q)$$
$$+ \frac{1}{q-1}\{R(m;q)R(m_1;q)R(m_2;q) - R(mm_1m_2;q)\}$$

if $q \nmid r$ because $H(0; q) = 1$. This formula, in spite of being quite compact, is not convenient for extension to composite moduli since it lacks multiplicative properties.

Let $q$ be squarefree. Notice that (10.4) and (10.5) are purely multiplicative in $q$ while (10.6) satisfies the twisted multiplication rule

(10.11) $$H(w; q_1 q_2) = H(w\bar{q}_1; q_2) H(w\bar{q}_2; q_1)$$

if $q = q_1 q_2$ (here $H(\star; q)$ is defined by (10.7) with $\chi = \chi_q$ for $q = q_1, q_2$ and $q_1 q_2$ respectively). Put

(10.12) $$h = (r, q), \ k = (mm_1 m_2, q/h), \ \ell = q/hk.$$

We write

$$H_r(m, m_1, m_2; q) = H_r(\overline{k\ell}m, m_1, m_2; h) H_r(\overline{\ell h}m, m_1, m_2; k) H_r(\overline{hk}m, m_1, m_2; \ell),$$

by the rule (10.3). Applying Lemma 10.1 we arrive at

(10.13)
$$H_r(m, m_1, m_2; q) = \chi_h(-1) \frac{h}{\phi(k)} R(m; k) R(m_1; k) R(m_2; k) H(\overline{rhk}mm_1 m_2; \ell)$$

provided $(h, mm_1 m_2) = 1$, or else $H_r(m, m_1, m_2; q)$ vanishes. Combining (10.13) and (10.1) we obtain:

LEMMA 10.2. *Let $c = qr$ with $q$ squarefree. Suppose $m, m_1, m_2$ are integers with*

(10.14) $$(m, r) = 1, \ (m_1 m_2, q, r) = 1.$$

*Then the modified character sum (8.34) satisfies*

(10.15)
$$G'(m, m_1, m_2; c) = \chi_{k\ell}(-1) \frac{r^2 qh}{\phi(k)} R(m; k) R(m_1; k) R(m_2; k) H(\overline{rhk}mm_1 m_2; \ell)$$

*where $h = (r, q), k = (mm_1 m_2, q)$ and $\ell = q/hk$. If the co-primality conditions (10.14) are not satisfied then $G'(m, m_1, m_2; c)$ vanishes.*

Applying trivial estimates $|R(m; k)| \leq \phi(k)$ and $|H(w; \ell)| \leq \ell^2$ we obtain:

COROLLARY 10.3. *Let $c = qr$ with $q$ squarefree. For any $m, m_1, m_2$ we have*

(10.16) $$|G(m, m_1, m_2; c)| \leq q^3 r^2.$$

## 11. Bilinear forms $\mathcal{H}$

Recall that $H(w; q)$ is the character sum defined by (10.7) where $\chi = \chi_q$ is the real primitive character of conductor $q$ (the Jacobi symbol). In this section we estimate general sums of type

$$(11.1) \qquad \mathcal{H} = \sum\sum\sum_{(rmn,q)=1} \gamma_{r,m} \beta_n H(a\bar{r}mn; q).$$

with any complex coefficients $\gamma_{r,m}$ and $\beta_n$ for $1 \leqslant r \leqslant R, 1 \leqslant m \leqslant M$ and $1 \leqslant n \leqslant N$. Although (11.1) is a triple sum we consider it as a bilinear form since the variables $r$ and $m$ are not separated. For convenience we assume that

$$(11.2) \qquad |\gamma_{r,m}| \leqslant 1,$$

but we make no conditions about $\beta = (\beta_n)$. We shall estimate $\mathcal{H}$ in terms of the $\ell_2$-norm

$$(11.3) \qquad \|\beta\| = (\sum_n |\beta_n|^2)^{\frac{1}{2}}.$$

LEMMA 11.1. *Let $(a, q) = 1$. Now,*

$$(11.4) \qquad \mathcal{H} \ll \|\beta\| (q + RM)^{\frac{1}{2}} (q + N)^{\frac{1}{2}} (qRM)^{\frac{1}{2}+\varepsilon}$$

*for any $\varepsilon > 0$, the implied constant depending only on $\varepsilon$.*

Actually we first prove (11.4) for the following sum

$$(11.5) \qquad \mathcal{H}^* = \sum\sum\sum_{(rmn,q)=1} \gamma_{r,m} \beta_n H^*(a\bar{r}mn; q)$$

where $H^*(w; q)$ denotes the reduced character sum

$$(11.6) \qquad H^*(w; q) = \sum\sum_{\substack{u,v \,(\mathrm{mod}\, q) \\ (uv-1,q)=1}} \chi(uv(u+1)(v+1)) e_q((uv-1)w).$$

Let $H(q)$ denote the above character sum restricted by the condition $uv \equiv 1(\mathrm{mod}\ q)$ instead of $(uv-1, q) = 1$. For $q$ prime $H(q) = -\chi(-1)$, whence for any squarefree $q$ we obtain $H(q) = \mu(q)\chi_q(-1)$ by the pure multiplicativity. Then by the twisted multiplicativity we derive

$$(11.7) \qquad H(w; q) = \sum\sum_{q_1 q_2 = q} \mu(q_1) \chi_{q_1}(-1) H^*(\bar{q}_1 w; q_2).$$

Hence it is clear that (11.4) for $\mathcal{H}^*$ implies that for $\mathcal{H}$.

Now we can express the additive character $e_q((uv-1)w)$ in $\mathcal{H}^*$ by means of the multiplicative characters $\psi(\mathrm{mod}\ q)$. Indeed for $(a, q) = 1$,

$$(11.8) \qquad e_q(a) = \frac{1}{\phi(q)} \sum_{\psi (\mathrm{mod}\, q)} \tau(\psi) \bar\psi(a)$$

where $\tau(\psi)$ is the Gauss sum. Hence we obtain

$$(11.9) \qquad H^*(w; q) = \frac{1}{\phi(q)} \sum_{\psi (\mathrm{mod}\, q)} \tau(\bar\psi) g(\chi, \psi) \psi(w)$$

where $g(\chi, \psi)$ is the hybrid character sum

$$(11.10) \qquad g(\chi, \psi) = \sum\sum_{u,v (\mathrm{mod}\, q)} \chi(uv(u+1)(v+1)) \psi(uv - 1).$$

Inserting (11.9) into (11.5) we obtain
$$(11.11)$$
$$\mathcal{H}^* = \frac{1}{\phi(q)} \sum_{\psi (\mathrm{mod}\, q)} \tau(\bar\psi) g(\chi, \psi) \psi(a) \left( \sum_r \sum_m \gamma_{r,m} \psi(\bar r m) \right) \left( \sum_n \beta_n \psi(n) \right).$$

In the last two sections we shall prove that for any $\psi (\mathrm{mod}\, q)$

$$(11.12) \qquad g(\chi, \psi) \ll q^{1+\varepsilon}$$

where the implied constant depends only on $\varepsilon$. Since $|\tau(\psi)| \leqslant q^{\frac{1}{2}}$ we obtain

$$|\mathcal{H}^*| \ll q^{\frac{1}{2}+\varepsilon} \sum_{\psi (\mathrm{mod}\, q)} |\sum_r \sum_m \gamma_{r,m} \psi(\bar r m)| |\sum_n \beta_n \psi(n)|.$$

Finally by Cauchy's inequality and the orthogonality of characters this yields (11.4) for $\mathcal{H}^*$. This completes the proof of Lemma 11.1.

*Remarks.* The bound (11.12), which is best possible, is essential. We need it for all $\psi (\mathrm{mod}\, q)$; even one exception would weaken the final results considerably. For example, if for one character we only had $g(\chi, \psi) \ll q^{\frac{3}{2}+\varepsilon}$, then the method would yield $L_f(\frac{1}{2}, \chi) \ll q^{\frac{3}{8}+\varepsilon}$ in place of (1.6), and $L(s, \chi) \ll q^{\frac{3}{16}+\varepsilon}$ in place of (1.12). Knowing that $H(w; q) \ll q$ (by the Riemann hypothesis for varieties) one can see that the estimate for bilinear form (11.4) saves an extra factor $q^{\frac{1}{2}}$ because of cancellation which is due to a variation in the angle of $H(w; q)$ as the parameter $w$ ranges over special numbers. Moreover it is very important that the three variables $r, m, n$ appear only in one block $w = a\bar r m n$ in a multiplicative fashion. The point is that we lost a factor $q^{\frac{1}{2}}$ only one time when passing from additive to multiplicative characters while we gained this factor twice when applying the orthogonality of the multiplicative characters. Another interesting point is that $r$ originated from moduli of Kloosterman sums $S(m, n; qr)$ and was transformed to a variable modulo $q$ by means of a kind of reciprocity.

## 12. Estimation of $\mathcal{G}(M, M_1, M_2; C)$

Inserting (10.15) into (8.36) we obtain

$$\mathcal{G}(M, M_1, M_2; C) \leqslant \frac{1}{q} \sum_{hk\ell=q} \sum \sum \frac{h}{\phi(k)} \sum_{\substack{(r, k\ell)=1 \\ C \leqslant hqr < 2C}}$$

$$\cdot \Big| \sum_{\substack{(mm_1m_2, h\ell)=1 \\ mm_1m_2 \equiv 0(k) \\ (m,r)=1}} \sum \sum R(m; k) R(m_1; k) R(m_2; k) H(\overline{rh^2k}mm_1m_2; \ell) \alpha_m \beta_{m_1, m_2} \Big|.$$

We put $m$ into outer summation, estimate $R(m, k)$ by $(m, k)\phi(k)/k$, and keep $m_1, m_2$ in the inner summation. To simplify the inner sum we create one variable $n = m_1 m_2$ with coefficients

(12.1) $$\beta_n(k) = \sum_{m_1 m_2 = n} \beta_{m_1, m_2} R(m_1; k) R(m_2; k).$$

We obtain

(12.2) $$\mathcal{G}(M, M_1, M_2; C) \leqslant \frac{1}{q} \sum_{hk\ell=q} \sum \sum \frac{h}{k} \sum_{\substack{(r, k\ell)=1 \\ C \leqslant hqr < 2C}} \sum_{\substack{|m| \leqslant M \\ (m, h\ell r)=1}} (m, k)$$

$$\cdot \Big| \sum_{\substack{(n, h\ell)=1 \\ mn \equiv 0(k)}} H(\overline{rh^2k}mn; \ell) \beta_n(k) \Big|.$$

Before making further simplifications we isolate the contribution to $\mathcal{G}(M, M_1, M_2; C)$ of the terms with $mm_1m_2 = 0$; we denote this contribution by $\mathcal{G}_0(M, M_1, M_2; C)$. For these terms on the right side of (12.2) we have $h = \ell = 1$ and $k = q$; therefore

$$\mathcal{G}_0(M, M_1, M_2; C) \ll Cq^{-3} \sum \sum \sum_{mm_1m_2=0} (m, q)(m_1, q)(m_2, q)$$

where the summation is also restricted by $|m| \leqslant M$, $|m_1| \leqslant M_1$, $|m_2| \leqslant M_2$. Hence

(12.3) $$\mathcal{G}_0(M, M_1, M_2; C) \ll C \left(1 + (M + M_1 + M_2)q^{-1}\right)^2 \tau^2(q).$$

Let $\mathcal{G}^*(M, M_1, M_2; C)$ denote the contribution to $\mathcal{G}(M, M_1, M_2; C)$ of the terms with $mm_1m_2$ different from zero. Putting $\kappa = (m, k)$ in (12.2) we get

$$|\mathcal{G}^*(M, M_1, M_2; C)| \leqslant \frac{1}{q} \sum_{hk\ell=q} \frac{h}{k} \sum_{\kappa\nu=k} \kappa \sum_{\substack{R \leqslant r < 2R \\ (r, \ell)=1}} \sum_{\substack{0 < |m| \leqslant M' \\ (m, \ell)=1}}$$

$$\cdot \Big| \sum_{\substack{0 < |n| \leqslant N' \\ (n, h\ell)=1}} H(\overline{rh^2}mn; \ell) \beta_{\nu n}(k) \Big|$$

where $R = C/hq, M' = M/\kappa, N' = M_1 M_2/\nu$ and $\beta_{\nu n}(k)$ is as in (12.1). For $n \neq 0$ we have $|\beta_n(k)| \leqslant (n,k)^2 \tau(n)$, so that

$$\sum_{0<|n|\leqslant N'} |\beta_{\nu n}(k)|^2 \leqslant \nu^4 \tau^2(\nu) \sum_{0<|n|\leqslant N'} (n,\kappa)^4 \tau^2(n)$$

$$\ll \nu^4 \tau^2(\nu) \kappa^3 N' (\log 2N')^4 \ll k^3 \tau^2(k) M_1 M_2 (\log 2 M_1 M_2)^4.$$

Now applying (11.4) we derive that
(12.4)
$$\mathcal{G}^*(M, M_1, M_2; C) \ll \left(1 + CMq^{-2}\right)^{\frac{1}{2}} \left(1 + M_1 M_2 q^{-1}\right)^{\frac{1}{2}} (CMM_1M_2)^{\frac{1}{2}+\varepsilon}.$$

Adding the estimates (12.3) and (12.4) we complete the proof of Proposition 8.2.

## 13. Estimation of $g(\chi, \psi)$

Our aim is to prove the bound (11.12) for the hybrid character sum $g(\chi, \psi)$ with $\chi(\bmod q)$ real, primitive character and $\psi(\bmod q)$ any character. Since $g(\chi, \psi)$ is multiplicative in the modulus it suffices to show:

LEMMA 13.1. *Let $p$ be prime, $\chi(\bmod p)$ the real, nonprincipal character, and $\psi(\bmod p)$ any character. Then*

(13.1) $$|g(\chi, \psi)| \ll p$$

*where the implied constant is absolute.*

In the proof we employ the Riemann hypothesis for varieties over the finite field $\mathbb{F}_p$. We need to consider the character sums over the field extensions $\mathbb{F}_q/\mathbb{F}_p$ with $q = p^m, m \geqslant 1$. Put

(13.2) $$g(\chi_m, \psi_m) = \sum\sum_{u,v\in\mathbb{F}_{p^m}} \chi_m(uv(u+1)(v+1))\psi_m(uv-1)$$

where $\chi_m$ and $\psi_m$ are the characters derived from $\chi, \psi$ by composing with the norm $N_{\mathbb{F}_q/\mathbb{F}_p} : \mathbb{F}_q \to \mathbb{F}_p$. Let $L(T)$ be the associated $L$-function

(13.3) $$L(T) = \exp\left(\sum_1^\infty g(\chi_m, \psi_m) \frac{T^m}{m}\right).$$

B. Dwork [Dw] showed that $L(T)$ is rational:

(13.4) $$L(T) = \prod_\nu (1 - \alpha_\nu T) \prod_\nu (1 - \beta_\nu T)^{-1}.$$

Equivalently

(13.5) $$g(\chi_m, \psi_m) = -\sum_\nu \alpha_\nu^m + \sum_\nu \beta_\nu^m.$$

P. Deligne [De] showed that the roots $\alpha_\nu^{-1}, \beta_\nu^{-1}$ are algebraic numbers with

(13.6) $$|\alpha_\nu| = p^{k_\nu/2}, \ |\beta_\nu| = p^{\ell_\nu/2}$$

where $k_\nu, \ell_\nu$ are nonnegative integers called the weights of $\alpha_\nu, \beta_\nu$. The total degree of $L(T)$ (the number of roots) was estimated by A. Adolphson and S. Sperber [AS] using Bombieri's idea [Bo]; it is bounded by a number independent of the characteristic $p$.

LEMMA 13.2. *Let $q = p^m$ and $\chi \in \hat{\mathbb{F}}_q, \chi \neq 1$. Then*

(13.7) $$\frac{1}{q-1} \sum_{\psi \in \hat{\mathbb{F}}_q} |g(\chi, \psi)|^2 = q^2 - 2q - 2.$$

*Proof.* By the orthogonality of characters the left-hand side of (13.7) is equal to

$$\sum\sum\sum_{u_1 v_1 = u_2 v_2}\sum \chi(u_1 v_1 (u_1+1)(v_1+1)) \bar\chi(u_2 v_2 (u_2+1)(v_2+1))$$

$$= \sum\sum\sum_{u_1, v_1, u_2 \in \mathbb{F}_q^*} \chi((u_1+1)(v_1+1))\bar\chi((u_2+1)(u_1 v_1 + u_2))\chi(u_2).$$

The sum over $u_1$ equals

$$\sum_{u_1 \neq 0} \chi(u_1+1)\bar\chi(u_1 v_1 + u_2) = \begin{cases} \bar\chi(v_1)(q-2) & \text{if } u_2 = v_1 \\ -\bar\chi(v_1) - \bar\chi(u_2) & \text{if } u_2 \neq v_1. \end{cases}$$

Now the sum over $u_2$ equals

$$\sum_{u_2 \neq 0} \bar\chi(u_2+1)\chi(u_2)\left(\sum_{u_1 \neq 0}\right) = q\bar\chi(v_1+1) + \bar\chi(v_1) + 1.$$

Finally the sum over $v_1$ equals

$$\sum_{v_1 \neq 0} \chi(v_1+1)\left(\sum_{u_2 \neq 0}\sum_{u_1 \neq 0}\right) = q(q-2) - 2.$$

This completes the proof of 13.7. □

Lemma 13.2 and the formula (13.5) imply (by choosing a suitable $m$) that all the weights $k_\nu, \ell_\nu$ are $\leq 2$ except for at most one root of weight three for one character $\psi$. If such a root exists then $\psi$ must be real, precisely $\psi = \chi$, because for the trivial character $\psi_0 = 1$ we have $g(\chi, \psi_0) = 1$. This completes the proof of Lemma 13.1 for all $\psi \neq \chi$. For $\psi = \chi$ we shall estimate $g(\chi, \chi)$ without recourse to the Riemann hypothesis.

## 14. Estimation of $g(\chi, \chi)$

There are several interesting expressions for $g(\chi, \chi)$ with real character $\chi \pmod p$. First recall that by the definition

$$(14.1) \qquad g(\chi, \chi) = \sum\sum_{u,v \pmod p} \chi(uv(u+1)(v+1)(uv-1)).$$

Changing variables one gets

$$(14.2) \qquad g(\chi, \chi) = \sum\sum_{u,v \pmod p} \chi(uv(u+1)(v+1)(u-v)).$$

Therefore $g(\chi, \chi)$ is associated with the elliptic curves in the Legendre family

$$(14.3) \qquad E_\lambda : y^2 = x(x+1)(x+\lambda).$$

Next we write (14.2) in the following way:

$$(14.4) \qquad g(\chi, \chi) = \sum\sum_{u,v \pmod p} \chi(2(u^2-1)(v^2-1)(u-v)).$$

Notice that $g(\chi, \chi) = \chi(-1)g(\chi, \chi)$, i.e. $g(\chi, \chi) = 0$ if $\chi(-1) = -1$. Now we separate the variables in (14.4) by inserting

$$(14.5) \qquad \chi(2(u-v)) = \frac{1}{\tau(\chi)} \sum_{a \pmod p} \chi(a)e_p(2a(u-v)).$$

This yields

$$(14.6) \qquad g(\chi, \chi) = \frac{1}{\tau(\chi)} \sum_a \chi(a) \left| \sum_u \chi(u^2-1)e_p(2au) \right|^2.$$

The innermost sum is a Kloosterman sum, precisely

$$(14.7) \qquad \sum_u \chi(u^2-1)e_p(2au) = S(a,a;p)$$

if $a \not\equiv 0 \pmod p$. Therefore we have

$$(14.8) \qquad g(\chi, \chi)\tau(\chi) = \sum_a \chi(a)S^2(a,a;p).$$

Weil's estimate for the Kloosterman sums $S(a, a; p)$, or Hasse's estimate for the number of points on the curve $E_\lambda/\mathbb{F}_p$, yields $|g(\chi, \chi)| \leqslant 2p^{3/2}$ while our goal is:

LEMMA 14.1. *For the real character $\chi(\bmod p), \chi \neq 1$,*

(14.9) $$|g(\chi, \chi)| \leqslant 4p.$$

Writing $S(a, a; p) = 2\sqrt{p}\, \cos(\omega_p(a))$ we derive from (14.8) that

$$g(\chi, \chi) = \tau(\chi)(1 + \chi(-1)) \sum_{a(\bmod p)} \chi(a) e^{2i\omega_p(a)}.$$

Therefore the estimate (14.9) reads as

(14.10) $$\Big| \sum_{a(\bmod p)} \chi(a) e^{2i\omega_p(a)} \Big| \leqslant 2\sqrt{p}.$$

In other words one can say that the variation of the Kloosterman angle $\omega_p(a)$ is quite independent of the sign change of the character $\chi(a)$.

Our proof of (14.9) is completely elementary. First notice that $S(a, a; p) = S(ab, a\bar{b}; p)$ for any $b \not\equiv 0 (\bmod p)$. By means of this parameter we create $p - 1$ copies of $g(\chi, \chi)$. Moreover changing $a$ into $ac^2$ with $c \not\equiv 0 (\bmod p)$ we create altogether $(p - 1)^2$ copies of $g(\chi, \chi)$. Thus,

$$p(p-1)^2 g^2(\chi, \chi) = {\sum_b}^\star {\sum_c}^\star \Big( \sum_a \chi(a) S^2(abc^2, a\bar{b}c^2; p) \Big)^2.$$

Put $bc^2 = x$ and $\bar{b}c^2 = y$, i.e. $b^2 = x/y$ and $c^4 = xy$. Given $x$ and $y$ there are at most two solutions for $b$ and four for $c$; therefore $g^2(\chi, \chi) \leqslant T$, where

(14.11) $$T = \frac{8}{p(p-1)^2} \sum_x \sum_y \Big( \sum_a \chi(a) S^2(ax, ay; p) \Big)^2.$$

We shall compute $T$ exactly. Opening the Kloosterman sum, squaring out and executing the summation in $x, y (\bmod p)$ we get

(14.12) $$T = \frac{8p}{(p-1)^2} \sum_{a_1} \sum_{a_2} \chi(a_1 a_2) \nu(a_1/a_2)$$

where $\nu(a)$ is the number of solutions to the system

$$a(d_1 + d_2) = d_3 + d_4, \qquad a\left(\frac{1}{d_1} + \frac{1}{d_2}\right) = \frac{1}{d_3} + \frac{1}{d_4}.$$

If $d_1 + d_2 = 0$ then $d_3 + d_4 = 0$. The number of such solutions is $(p-1)^2$. If $d_1 + d_2 \neq 0$ then the second equation of the system can be replaced by $d_1 d_2 = d_3 d_4$. Given $d_1, d_2$ with $d_1 d_2 (d_1 + d_2) \neq 0$ the $d_3, d_4$ are the roots of $X^2 - a(d_1 + d_2)X + d_1d_2 = 0$, and the number of roots equals $1 + \chi[a^2(d_1 + d_2)^2 - 4d_1d_2]$. Hence

$$\nu(a) = (p-1)^2 + \sum\sum_{d_1d_2(d_1+d_2)\neq 0} (1 + \chi[a^2(d_1+d_2)^2 - 4d_1d_2])$$

$$= (p-1)^2 + (p-1) \sum_{d(d+1)\neq 0} (1 + \chi[a^2(d+1)^2 - 4d])$$

$$= (p-1)^2 + (p-1)(p-2) + (p-1)\sum_{d(d-1)\neq 0} \chi(a^2d^2 - 4d + 4).$$

The last sum is equal to $-3 + p\delta(a^2)$, where $\delta(a) = 1$ if $a = 1$, and $\delta(a) = 0$ otherwise, therefore

(14.13) $$\nu(a) = 2(p-1)(p-3) + p(p-1)\delta(a^2).$$

Inserting (14.13) into (14.12) we arrive at

(14.14) $$T = 8p^2(1 + \chi(-1)).$$

This, together with $g^2(\chi, \chi) \leqslant T$, completes the proof of (14.9).

*Remarks.* Some of the ideas of our proof of Lemma 14.1 are reminiscent of the Kloosterman arguments [Kl]. The method is capable of producing good results for the sum $\sum \chi(a)S^2(a, a; p)$ with any character $\chi$ which assumes one value with large multiplicity. In particular it works well for a character of any fixed order. We believe that (14.10) is true for any $\chi(\mathrm{mod}\, q)$. A more advanced study of the Kloosterman angles $\omega_p(a)$ can be found in the book [Ka].

*Added in proof.* W. Duke showed (in February, 1999) that $g(\chi, \chi) = 2\mathrm{Re}J^2(\chi, \psi)$ where $J(\chi, \psi)$ is the Jacobi sum and $\psi$ is a quartic character modulo $p \equiv 1 (\mathrm{mod}\, 4)$. He also pointed out that $g(\chi, \chi)$ is the $p^{\mathrm{th}}$ Fourier coefficient of $\eta(4z)^6 \in S_3(\Gamma_0(12))$, so Lemma 14.1 follows from Ramanujan's conjecture (proved by P. Deligne). N. Katz also informed us that $g(\chi, \chi)$ can be computed explicitly.

AMERICAN INSTITUTE OF MATHEMATICS, PALO ALTO, CA
OKLAHOMA STATE UNIVERSITY, STILLWATER, OK
*E-mail address*: conrey@aimath.org

RUTGERS UNIVERSITY, PISCATAWAY, NJ
*E-mail address*: iwaniec@math.rutgers.edu

References

[AS]   A. ADOLPHSON and S. SPERBER, Character sums in finite fields, *Compositio Math.* **52** (1984), 325–354.
[AL]   A. ATKIN and J. LEHNER, Hecke operators on $\Gamma_0(m)$, *Math. Ann.* **185** (1970), 134–160.
[Bo]   E. BOMBIERI, On exponential sums in finite fields, II, *Invent. Math.* **47** (1978), 29–39.
[Bu]   D. A. BURGESS, On character sums and $L$-series, I, *Proc. London Math. Soc.* **12** (1962), 193–206.
[By]   V. A. BYKOVKY, Estimation of Hecke series for holomorphic parabolic forms twisted by a character on the critical line (in Russian), Preprint 05-1995 of the Institute for Applied Mathematics of the Russian Academy of Sciences in Khabarovsk.
[Da]   H. DAVENPORT, *Multiplicative Number Theory*, Markham Publ. Co., Chicago, 1967.
[De]   P. DELIGNE, La conjecture de Weil, II, *Publ. Math. I.H.E.S.* **52** (1981), 137–252.
[DI]   J.-M. DESHOUILLERS and H. IWANIEC, Kloosterman sums and Fourier coefficients of cusp forms, *Invent. Math.* **70** (1982), 219–288.
[DFI]  W. DUKE, J. B. FRIEDLANDER, and H. IWANIEC, Bounds for automorphic $L$-functions, *Invent. Math.* **112** (1993), 1–8.
[Dw]   B. DWORK, On the rationality of the zeta function of an algebraic variety, *Amer. J. Math.* **82** (1960), 631–648.
[Fr]   J. B. FRIEDLANDER, Bounds for $L$-functions, *Proc. I.C.M.* (Zürich, 1994), Birkhäuser, Basel, 1995, 363–373.
[Gu]   J. GUO, On the positivity of the central critical values of automorphic $L$-functions for GL(2), *Duke Math. J.* **83** (1996), 157–190.
[H-B]  D. R. HEATH-BROWN, Hybrid bounds for Dirichlet $L$-function, II, *Quart. J. Math. Oxford* **31** (1980), 157–167.
[He]   D. A. HEJHAL, *The Selberg Trace Formula for* PSL$(2, \mathbb{R})$, *Lecture Notes in Math.* **1001** (1983), Springer-Verlag, New York.
[HL]   J. HOFFSTEIN and P. LOCKHART, Coefficients of Maass forms and the Siegel zero, *Ann. of Math.* **140** (1994), 161–181.
[Hu]   M. N. HUXLEY, Exponential sums and the Riemann zeta function, IV, *Proc. London Math. Soc.* **66** (1993), 1–40.
[I1]   H. IWANIEC, Fourier coefficients of modular forms of half-integral weight, *Invent. Math.* **87** (1987), 385–401.
[I2]   ———, *Introduction to the Spectral Theory of Automorphic Forms*, Biblioteca de la Revista Matemática Iberoamericana (1995), Madrid.
[I3]   ———, *Topics in Classical Automorphic Forms*, *Graduate Studies in Math.* **17**, Providence, RI (1997).
[I4]   ———, Small eigenvalues of Lapacian for $\Gamma_0(N)$, *Acta Arith.* **56** (1990), 65–82.
[ILS]  H. IWANIEC, W. LUO, and P. SARNAK, Low lying zeros of families of $L$-functions, *Publ. Math. I.H.É.S*, to appear.
[IS]   H. IWANIEC and P. SARNAK, The non-vanishing of central values of automorphic $L$-functions and Siegel's zero, preprint, 1997.
[KS]   S. KATOK and P. SARNAK, Heegner points, cycles, and Maass forms, *Israel J. Math.* **84** (1993), 193–227.
[Ka]   N. M. KATZ, *Gauss Sums, Kloosterman Sums, and Monodromy Groups*, *Annals of Math. Studies* **116** (1988), Princeton University Press, Princeton, NJ.
[Kl]   H. D. KLOOSTERMAN, On the representation of numbers in the form $ax^2 + by^2 + cz^2 + dt^2$, *Acta Math.* **49** (1926), 407–464.
[KZ]   W. KOHNEN and D. ZAGIER, Values of $L$-series of modular forms at the center of the critical strip, *Invent. Math.* **64** (1981), 175–198.
[L]    W. LI, Newforms and functional equations, *Math. Ann.* **212** (1975), 285–315.
[R]    MICHAEL J. RAZAR, Modular forms for $G_0(N)$ and Dirichlet series, *Trans. A.M.S.* **231** (1977), 489–495.

[Wa]   J.-L. Waldspurger, Sur les coefficients de Fourier des forms modulaires de poids demi-entier, *J. Math. Pures Appl.* **60** (1981), 375–484.
[W]    G. N. Watson, *A Treatise on the Theory of Bessel Functions*, Cambridge University Press, Cambridge, England (edition 1995).